\documentclass[11pt]{amsart}

\usepackage{amssymb}
\usepackage[all]{xy}
\xyoption{ps}
\xyoption{dvips}
\usepackage{amsthm}

\newtheorem{Thm}{Theorem}[section]
\newtheorem{Lem}[Thm]{Lemma}
\newtheorem{Cor}[Thm]{Corollary}
\newtheorem{Prop}[Thm]{Proposition}

\newcommand{\rk}{\operatorname{rk}}
\newcommand{\md}{\operatorname{mod}}

\newcommand{\Hom}{\operatorname{Hom}}
\newcommand{\Ext}{\operatorname{Ext}}
\newcommand{\End}{\operatorname{End}}

\newcommand{\dm}{{\rm dim}\,}
\newcommand{\M}{{\rm M}}
\newcommand{\GL}{{\rm GL}}
\newcommand{\field}{K}
\newcommand{\orb}{{\mathcal O}}
\newcommand{\V}{{\rm V}}
\newcommand{\sub}{{\mathcal S}}
\newcommand{\imod}{{\mathcal I}}
\newcommand{\nil}{{\rm N}}

\newcommand{\ke}{\operatorname{Ker}}
\newcommand{\im}{\operatorname{Im}}
\newcommand{\soc}{\operatorname{soc}}
\newcommand{\tp}{\operatorname{top}}

\newcommand{\reg}{{\rm reg}}

\begin{document}

\title[Pairs of nilpotent matrices annihilating each other]
{Varieties of pairs of nilpotent matrices annihilating each other}
\author{Jan Schr\"oer}
\address{Jan Schr\"oer\newline
Department of Pure Mathematics\newline
University of Leeds\newline
Leeds LS2 9JT\newline
ENGLAND}
\email{jschroer@maths.leeds.ac.uk}


\begin{abstract}
We classify the irreducible components of the varieties
\[ 
\V(n,a,b) = \{ (A,B) \in \M_n(\field) \times \M_n(\field) \mid 
AB = BA = A^a = B^b = 0 \}. 
\]
\end{abstract}

\maketitle
\tableofcontents

\section{Introduction and main results}


Let $\M_n(\field)$ be the set of $n \times n$-matrices with
entries in an algebraically closed field $\field$.
The study of affine varieties given by matrices or pairs of matrices,
which satisfy certain relations, is a classical subject.
One fundamental question is the decomposition of these varieties into 
irreducible components. 
Consider the varieties 
\[ 
\nil(n,l) = \{ M \in \M_n(\field) \mid M^l = 0 \} 
\]
and
\[ 
{\rm Z}(n) = 
\{ (A,B) \in \M_n(\field) \times \M_n(\field) \mid AB = BA = 0 \}. 
\]
The variety $\nil(n,l)$ is irreducible by \cite{Ge} and \cite{H}, and
the irreducible components of ${\rm Z}(n)$ are
\[ 
\{ (A,B) \in {\rm Z}(n) \mid \rk(A) \leq n-i, \rk(B) \leq i \} 
\]
for $0 \leq i \leq n$.
For $n,a,b \geq 2$ define
\begin{align*}
\V(n,a,b) &= 
\{ (A,B) \in \M_n(\field) \times \M_n(\field) \mid 
AB = BA = A^a = B^b = 0 \}\\
&= \left( \nil(n,a) \times \nil(n,b) \right) \cap {\rm Z}(n). 
\end{align*}
Our main result is the classification of irreducible components of
$\V(n,a,b)$.
This question appears for $a = b = n$ as an open problem in 
\cite[Problem 3, p.208]{K}. 
In this special case, we get the following surprising result:

\begin{Thm}\label{nil2}
The irreducible components of $\V(n,n,n)$ are
\[ 
\{ (A,B) \in \V(n,n,n) \mid \rk(A) \leq n-i, \rk(B) 
\leq i \} 
\]
for $1 \leq i \leq n-1$.
Each component has dimension $n^2-n+1$.
\end{Thm}

Thus each irreducible component of $\V(n,n,n)$ is the
intersection of an irreducible component of ${\rm Z}(n)$ with 
$\nil(n,n) \times \nil(n,n)$. 
The case $a=b=2$ and $n$ arbitrary was studied in \cite{M}.

A {\it partition} of $n$ is a sequence ${\bf p} = (p_1, \cdots, p_t)$ of 
positive integers such that $\sum_{i=1}^t p_i =  n$ and 
$p_i \geq p_{i+1}$ for all $i$.
Let $l({\bf p}) = t$ be the {\it length} of ${\bf p}$.
The set of partitions ${\bf p}$ of $n$ with $p_i \leq a$
for all $i$ is denoted by ${\mathcal P}(n,a)$.

By $\trianglelefteq$ we denote the usual dominance order on 
${\mathcal P}(n,a)$, 
see Section \ref{sectionregcomp} for a definition.

The conjugacy classes of matrices in $\nil(n,a)$ are parametrized by
${\mathcal P}(n,a)$.
Namely, for a matrix $M \in \nil(n,a)$, let $J(M)$ be its
Jordan normal form, and set $p(M) = (p_1, \cdots, p_t)$
where the $p_i$ are the sizes of the Jordan blocks of $J(M)$,
ordered decreasingly.
Clearly, we have $p(M) \in {\mathcal P}(n,a)$.
For ${\bf p} \in {\mathcal P}(n,a)$ 
let 
\[
C({\bf p}) = \{ M \in \nil(n,a) \mid p(M) = {\bf p} \}
\]
be the corresponding conjugacy class
in $\nil(n,a)$.

There are two projection maps
\[
\xymatrix{
 & \V(n,a,b) \ar[dr]^{\pi_2} \ar[dl]_{\pi_1} \\
\nil(n,a) & & \nil(n,b)
}
\]
where
$\pi_1(A,B) = A$ and
$\pi_2(A,B) = B$.
For ${\bf a} \in {\mathcal P}(n,a)$ let
\[
\Delta({\bf a}) = \pi_1^{-1}(C({\bf a})).
\]
In general, $\Delta({\bf a})$ is not irreducible.
Only if $a=b=n$, these sets have nice properties:

\begin{Thm}\label{nil1}
For each ${\bf a} \in {\mathcal P}(n,n)$ the set
$\Delta({\bf a}) \subset \V(n,n,n)$ is locally closed and irreducible.
We have
\[ 
\Delta(1, \cdots, 1) \subset \overline{\Delta(2,1, \cdots, 1)}, 
\]
and if ${\bf a} \not= (1, \cdots, 1)$, then 
\[ 
\Delta({\bf a}) \subset \overline{\Delta({\bf b})} 
\] 
if and only if ${\bf a} \trianglelefteq {\bf b}$ and 
$l({\bf a}) = l({\bf b})$.
\end{Thm}

For the study of the general case, define the {\it standard 
stratification} of $\V(n,a,b)$ as follows:
Let 
\[
{\mathcal P}(n,a,b) = {\mathcal P}(n,a) \times {\mathcal P}(n,b).
\]
For $({\bf a},{\bf b}) \in {\mathcal P}(n,a,b)$
let
\[
\Delta({\bf a},{\bf b}) = \pi_1^{-1}(C({\bf a})) \cap 
\pi_2^{-1}(C({\bf b}))
\]
be the corresponding stratum of the standard stratification.
Unfortunately, these strata are in general not very well-behaved:
\begin{itemize}

\item A stratum might be empty;

\item Strata are not necessarily irreducible;

\item The closure of a stratum is in general not a union of strata.

\end{itemize}
However, the socalled `regular strata' have nice properties.
Observe that for $(A,B) \in \V(n,a,b)$ the inequality 
\[ 
\rk(A) + \rk(B) \leq n 
\]
holds.
This follows already from the condition $AB = 0$.
We call $(A,B)$ {\it regular} if 
$\rk(A) + \rk(B) = n$.
An irreducible component of $\V(n,a,b)$ is {\it regular} if it contains
a regular element, and we
call $({\bf a},{\bf b}) \in {\mathcal P}(n,a,b)$ and 
also its corresponding stratum
$\Delta({\bf a},{\bf b})$ {\it regular} if $\Delta({\bf a},{\bf b})$
contains a regular element.

For a partition ${\bf p} = (p_1, \cdots, p_t) \not= (1, \cdots, 1)$ define
\[
{\bf p}-1 = (p_1-1, \cdots, p_r-1)
\]
where $r = \max \{ 1 \leq i \leq t \mid p_i \geq 2 \}$.
For example,
\[
(3,2,2,1) - 1 = (2,1,1).
\]
The following result determines which strata are regular.

\begin{Prop}\label{regularstrata}
For $({\bf a},{\bf b}) \in {\mathcal P}(n,a,b)$ the following are equivalent:
\begin{itemize}

\item[(1)]
$({\bf a},{\bf b})$ is regular;

\item[(2)]
$l({\bf a}) + l({\bf b}) = n$ and $l({\bf a}-1) = l({\bf b}-1)$.

\end{itemize}
In this case, all elements in $\Delta({\bf a},{\bf b})$ are regular.
\end{Prop}

If ${\bf p}$ is a partition, then let
\[
|i \in {\bf p}|
\]
be the number of entries of ${\bf p}$ which are equal to $i$.
The next theorem yields a classification of all
regular irreducible components.

\begin{Thm}\label{theorem1} 
If
$({\bf a},{\bf b}) \in {\mathcal P}(n,a,b)$ 
is regular, then
$\Delta({\bf a},{\bf b})$ is locally closed and irreducible.
In this case,
the closure of $\Delta({\bf a},{\bf b})$ is an irreducible component of 
$\V(n,a,b)$
if and only if the following hold:
\begin{itemize}

\item[(1)] ${\bf a}$ has at most one entry different from $1,2$ and $a$;

\item[(2)] ${\bf b}$ has at most one entry different from $1,2$ and $b$;

\item[(3)] $l({\bf a}-1) \leq |a \in {\bf a}| + |b \in {\bf b}| + 1$.

\end{itemize}
\end{Thm}

Next, we determine when all irreducible components of $\V(n,a,b)$ are 
regular.

\begin{Prop}\label{regulardense}
The set of regular elements is dense in $\V(n,a,b)$ if and only
if $n \leq a+b-2$ or $n = a+b$.
\end{Prop}

The classification of the non-regular irreducible components of
$\V(n,a,b)$ is more complicated and needs more notation.
We state and prove the result in Section \ref{sectionnonregcomp}.

The paper is organized as follows:
In Section \ref{varieties} we repeat some basics
on varieties of modules. 
In particular, we recall Richmond's construction 
of a stratification of these varieties, which we will
use throughout.
We regard $\V(n,a,b)$ as a variety of modules over a Gelfand-Ponomarev
algebra, and we use module theory to classify the irreducible
components of $\V(n,a,b)$.
Section \ref{stringalgebras} is a collection of mostly known
results on Gelfand-Ponomarev algebras.
Richmond's stratification turns out to be finite for
$\V(n,a,b)$.
This is studied in Section \ref{indexsection}.
In Section \ref{sectionregcomp} we prove that all regular strata are
irreducible.
This is used in Section \ref{nilcase1} to prove
Theorem \ref{nil2}.
The classification of all regular
components of $\V(n,a,b)$ can be found in Section \ref{classregcomp}.
Theorem \ref{nil1} is proved at the end of Section \ref{classregcomp}.
The main result of Section \ref{sectionnonregcomp} is the classification
of all non-regular components of $\V(n,a,b)$.
Finally, some examples are given in Section \ref{examples}.

\vspace{0.5cm}
\noindent 
{\bf Acknowledgements.}\,
The author received a Postdoctoral Fellowship from the DAAD, 
Germany, for a stay at the UNAM in Mexico City, where most of 
this work was done.
He thanks Christof Gei{\ss} and Lutz Hille for helpful and
interesting discussions.


\section{Varieties of modules}\label{varieties}


Let $A$ be a finitely generated $\field$-algebra.
Fix a set $a_1, \cdots, a_N$ of generators of $A$.
By $\md(A,n)$ we denote the affine variety of
$A$-module structures on $\field^n$.
Each such $A$-module structure corresponds to
a $\field$-algebra homomorphism $A \to \M_n(\field)$,
or equivalently to a tuple $(M_1, \cdots, M_N)$ of
$n \times n$-matrices such that the $M_i$ satisfy the same relations
as the $a_i$.
The group $\GL_n(\field)$ acts by simultaneous conjugation on
$\md(A,n)$, and the orbits of this action are in 1-1
correspondence with the isomorphism classes of $n$-dimensional $A$-modules.
An orbit $\orb(X)$ of a module $X$
has dimension $n^2 - \dm \End_A(X)$.
If $\orb(X)$ is contained in the closure of an orbit 
$\orb(Y)$, then we write $Y \leq_{\rm deg} X$.
It is well known that $Y \leq_{\rm deg} X$ implies
$\dm \Hom_A(Y,M) \leq \dm \Hom_A(X,M)$ for all modules $M$,
see for example \cite{Bo}.
If 
\[ 
0 \longrightarrow X \longrightarrow Y \longrightarrow Z 
\longrightarrow 0 
\] 
is a short exact sequence,
then $Y \leq_{\rm deg} X \oplus Z$.
If there exists a module $Z$ and a short exact sequence
\[ 
0 \longrightarrow X \longrightarrow Y \oplus Z \longrightarrow Z 
\longrightarrow 0, 
\]
then it is proved in \cite{Rie} that $Y \leq_{\rm deg} X$. 
The converse is also true by \cite{Z}.
Short exact sequences of this form  are called {\it Riedtmann sequences}.
We call a module $X$ a {\it minimal degeneration} if there exists no
module $Y$ with $Y <_{\rm deg} X$.

Now, let $A$ be a finite-dimensional $\field$-algebra, and
let $\imod_A(n)$ be a set of representatives of isomorphism 
classes of
submodules of $A^n$ which have dimension $n(d-1)$ where $d = 
\dm(A)$.
The modules in $\imod_A(n)$ are called the {\it index modules}
of $A$.
For each $L \in \imod_A(n)$ let $\sub(L)$ be the set of
points $X \in \md(A,n)$ 
such that there exists a short 
exact sequence
\[ 
0 \longrightarrow L \longrightarrow A^n \longrightarrow X
\longrightarrow 0 
\] 
of $A$-modules.
Such a set $\sub(L)$ is called a {\it stratum}.
Note that $\md(A,n)$
is the disjoint union of the
$\sub(L)$ where $L$ runs through $\imod_A(n)$.
The following theorem can be found in \cite{R}.

\begin{Thm}[Richmond]\label{richmond}
Let $A$ be a finite-dimensional $\field$-algebra.
Then the following hold:
\begin{enumerate}

\item[{\rm (1)}]
For each $L \in \imod_A(n)$ the
stratum $\sub(L)$ is smooth, locally closed, irreducible
and has dimension 
\[ 
\dm \Hom_A(L,A^n) - \dm \End_A(L);
\]
\item[{\rm (2)}]
Let $L, M \in \imod_A(n)$. 
If $\sub(L)$ is contained in the closure of
$\sub(M)$, then $M \leq_{\rm deg} L$;

\item[{\rm (3)}]
Let $L, M \in \imod_A(n)$.
If $M \leq_{\rm deg} L$ and
\[ 
\dim \Hom_A(L,A) = \dm \Hom_A(M,A), 
\] 
then $\sub(L)$ is contained in the closure of $\sub(M)$.

\end{enumerate}
\end{Thm}

Unfortunately, the converse of the second part of this theorem is
usually wrong.
So it remains a difficult problem to decide when a stratum is contained
in the closure of another stratum.
Another problem is, that the set $\imod_A(n)$ is often infinite.
Following \cite{R} 
an algebra $A$ is called {\it subfinite} if $\imod_A(n)$ is finite
for all $n$.


\section{Gelfand-Ponomarev algebras}\label{stringalgebras}


We identify $\V(n,a,b)$ with the variety of $n$-dimensional 
modules over the algebra 
\[ 
\Lambda = \Lambda_{a,b} = \field[x,y]/(xy,x^a,y^b). 
\]
We call $\Lambda$ a {\it Gelfand-Ponomarev algebra}.

The group $\GL_n(\field)$ acts on $\V(n,a,b) = \md(\Lambda,n)$ by 
simultaneous conjugation, i.e. 
\[ 
g \cdot (A,B) = (gAg^{-1},gBg^{-1}). 
\]
The orbits of this action are in 1-1 correspondence with the
isomorphism classes of $n$-dimensional $\Lambda$-modules.
By $\orb(M)$ we denote the orbit of an element $M \in \V(n,a,b)$.

In the following we repeat Gelfand and Ponomarev's 
classification of indecomposable $\Lambda$-modules 
(by a `module' we always mean a finite-dimensional right module).
As a main reference we use \cite{GP}, but see also \cite{BR}.

A {\it string} of length $n \geq 1$ is a word 
$c_1 \cdots c_n$ with letters $c_i \in \{ x,y \}$ such that
no subword is of the form $x^a$ or $y^b$.
Additionally, we define a string $1$ of length 0.
Set $x^0 = y^0 = 1$.

The length of an arbitrary string $C$ is denoted by $|C|$.
Let $C = c_1 \cdots c_n$ and $D = d_1 \cdots d_m$ be strings of length 
at least one. 
If $CD = c_1 \cdots c_n d_1 \cdots d_m$ is a string, then we say that the 
concatenation of $C$ and $D$ is defined.
For an arbitrary string $C$ let $1C = C1 = C$.

For each string $C$ we construct a {\it string module} $M(C)$ over $\Lambda$ 
as follows:
First, assume that $n = |C| \geq 1$.
Fix a basis $\{ z_1, \cdots ,z_{n+1} \}$ of $M(C)$. 
Given an arrow $\alpha \in \{ x,y \}$ let 
\[
z_i \cdot \alpha = 
\begin{cases}
z_{i+1} & \text{if $\alpha = c_i = y$ and $1 \leq i \leq n$},\\
z_{i-1} & \text{if $\alpha = c_{i-1} = x$ and $2 \leq i \leq n+1$},\\
0 & \text{otherwise}.
\end{cases}
\]
For $C = 1$ let $S = M(C)$ be the one-dimensional
module with basis $\{ z_1 \}$ such that
$z_1 \cdot x = z_1 \cdot y = 0$.
This is the unique simple $\Lambda$-module.
The $z_i$ are called the {\it canonical basis vectors} of $M(C)$.

\begin{figure}[ht]
\unitlength1.0cm
\begin{picture}(5,3)
\put(2,2.5){$z_3$}
\put(1.9,2.4){\vector(-1,-1){0.6}}
\put(1.4,2.2){$x$}
\put(1,1.5){$z_2$}
\put(0.9,1.4){\vector(-1,-1){0.6}}
\put(0.4,1.2){$x$}
\put(0,0.5){$z_1$}
\put(2.5,2.4){\vector(1,-1){0.6}}
\put(2.9,2.2){$y$}
\put(3,1.5){$z_4$}
\put(4,2.5){$z_5$}
\put(3.9,2.4){\vector(-1,-1){0.6}}
\put(3.4,2.2){$x$}
\put(4.5,2.4){\vector(1,-1){0.6}}
\put(4.9,2.2){$y$}
\put(5,1.5){$z_6$}
\end{picture}
\caption{The string module $M(xxyxy)$}\label{fig1}
\end{figure}
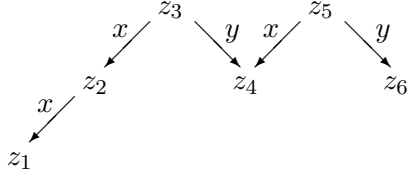

For example, let $C = xxyxy$.
Then the string module $M(C)$ looks as in Figure \ref{fig1},
where $z_1, \cdots, z_6$ are the canonical basis vectors of
$M(C)$, and the arrows indicate how the generators $x$ and $y$ of 
$\Lambda$ operate on these basis vectors. 
Set $(A,B) = M(xxyxy)$.
We have 
\[
(A,B) \in \pi_1^{-1}(C(3,2,1)) \cap \pi_2^{-1}(C(2,2,1,1)) = 
\Delta((3,2,1),(2,2,1,1)).
\]
A string $C = c_1 \cdots c_n$ of length at least one is called a {\it band}
if all powers $C^m$ are defined. 
Next, we associate to a given band $B = b_1 \cdots b_m$ and some 
$n \geq 1$ a family
\[ 
\{ M(B,\lambda_1, \cdots, \lambda_n) \mid \lambda_i \in \field^*,
1 \leq i \leq n \} 
\]
of {\it band modules}.
Fix a basis $\{ z_{1j}, \cdots , z_{mj} \mid 1 \leq j \leq n \}$
of $M(B,\lambda_1, \cdots, \lambda_n)$. 
For $\alpha \in \{ x,y \}$ define 
\[ 
z_{1j} \cdot \alpha = 
\begin{cases}
z_{2j} & \text{if $\alpha = b_1 = y$},\\
\lambda_j z_{mj} + z_{mj-1} & \text{if $\alpha = b_m = x$ and 
$2 \leq j \leq n$},\\
\lambda_1 z_{m1} & \text{if $\alpha = b_m = x$ and $j=1$},\\
0 & \text{otherwise},
\end{cases}
\]
and let
\[ 
z_{mj} \cdot \alpha = 
\begin{cases}
z_{m-1j} & \text{if $\alpha = b_{m-1} = x$},\\
\lambda_j z_{1j} + z_{1j-1} & \text{if $\alpha = b_m = y$ and 
$2 \leq j \leq n$},\\
\lambda_1 z_{11} & \text{if $\alpha = b_m = y$ and $j=1$},\\
0 & \text{otherwise}.
\end{cases} 
\]
For $2 \leq i \leq m-1$ and $1 \leq j \leq n$ we define
\[ 
z_{ij} \cdot \alpha = 
\begin{cases}
z_{i+1j} & \text{if $\alpha = b_i = y$},\\
z_{i-1j} & \text{if $\alpha = b_{i-1} = x$},\\
0 & \text{otherwise}.
\end{cases}  
\]
The $z_{ij}$ are called the {\it canonical basis vectors} of 
$M(B, \lambda_1, \cdots, \lambda_n)$.

\begin{figure}[ht]
\unitlength1.0cm
\begin{picture}(7,4.5)

\put(0.5,0.5){$z_{11}$} 
\put(5.3,0.6){\vector(-1,0){4.2}}\put(2.8,0.3){$(y,\lambda_1)$}
\put(5.7,0.8){\vector(0,1){3.0}}\put(5.8,2.2){$x$}
\put(5.5,0.5){$z_{51}$}

\put(1.5,1.5){$z_{12}$}
\put(4.3,1.6){\vector(-1,0){2.2}}\put(2.8,1.7){$(y,\lambda_2)$}
\put(4.7,1.8){\vector(0,1){1.0}}\put(4.8,2.2){$x$}
\put(4.5,1.5){$z_{52}$}
\put(4.3,1.5){\vector(-4,-1){3.2}}\put(3.2,1.0){$y$}

\put(1.5,3){$z_{22}$}
\put(1.7,2.8){\vector(0,-1){1.0}}\put(1.4,2.2){$x$}
\put(2.8,3.1){\vector(-1,0){0.7}}\put(2.5,3.2){$x$}
\put(3,3){$z_{32}$}
\put(3.6,3.1){\vector(1,0){0.7}}\put(3.8,3.2){$y$}
\put(4.5,3){$z_{42}$}

\put(0.5,4){$z_{21}$}
\put(0.7,3.8){\vector(0,-1){3.0}}\put(0.4,2.2){$x$}
\put(2.8,4.1){\vector(-1,0){1.7}}\put(2.0,4.2){$x$}
\put(3,4){$z_{31}$}
\put(3.6,4.1){\vector(1,0){1.7}}\put(4.1,4.2){$y$}
\put(5.5,4){$z_{41}$}

\end{picture}
\caption{The band module $M(xxyxy,\lambda_1,\lambda_2)$}\label{fig2}
\end{figure}
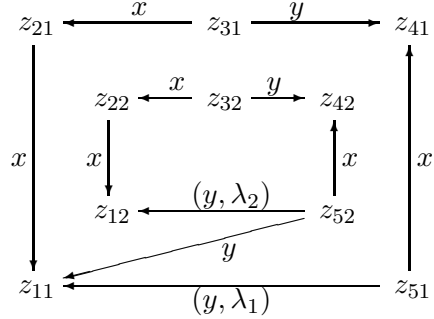

For example, let $B = xxyxy$.
Then the band module $M(B,\lambda_1,\lambda_2)$ looks as in Figure 
\ref{fig2}.
The arrows in Figure \ref{fig2} indicate how the generators $x$ and $y$ of
$\Lambda$ operate on the canonical basis vectors of 
$M(B,\lambda_1,\lambda_2)$.
For example, $z_{51} \cdot y = \lambda_1 z_{11}$,
$z_{52} \cdot y = \lambda_2 z_{12} + z_{11}$,
$z_{32} \cdot y = z_{42}$ etc.

The next lemma is proved by straightforward base change 
calculations.

\begin{Lem}
Let $M(B,\lambda_1, \cdots, \lambda_n)$ be a band module. 
If $\lambda_l \not= \lambda_{l+1}$ for some $l$, then 
$M(B,\lambda_1, \cdots, \lambda_n)$ is isomorphic to 
\[ 
M(B,\lambda_1, \cdots, \lambda_l) \oplus 
M(B,\lambda_{l+1}, \cdots, \lambda_n). 
\]
\end{Lem}

If $\lambda_i = \lambda_j$ for all $i$ and $j$, then define
$M(B,\lambda,n) = M(B,\lambda_1, \cdots, \lambda_n)$,
compare \cite{BR}.

A band $B$ is called {\it periodic} if there exists some 
string $C$ such that $B = C^m$ for some $m \geq 2$.
A band is called {\it primitive} if 
it is not periodic.
For primitive bands $B_1$ and $B_2$ define 
$B_1 \sim B_2$ if $B_1 = BB'$ and $B_2 = B'B$ for some strings
$B$ and $B'$.
Let $\sub$ be the set of strings, and let
${\mathcal B}$ be a set of representatives of 
equivalence classes of primitive bands with 
respect to the equivalence relation $\sim$.
The following theorem is proved in \cite{GP}.

\begin{Thm}[Gelfand-Ponomarev]\label{indecomposables}
The modules $M(C)$ and $M(B,\lambda,n)$ with
$C \in \sub$, $B \in {\mathcal B}$, $\lambda \in \field^*$ and 
$n \geq 1$ is a complete set of representatives of isomorphism classes
of indecomposable $\Lambda$-modules.
\end{Thm}

The next lemma follows from the construction of string and
band modules and from Theorem \ref{indecomposables}.

\begin{Lem}\label{rankofrepresentations}
If $(A,B) \in \V(n,a,b)$, then 
\[ 
n - s = \rk(A) + \rk(B), 
\]
where $s$ is the number of string modules in a decomposition of $(A,B)$
into a direct sum of indecomposable modules.
\end{Lem}

\begin{Cor}
An element in $\V(n,a,b)$ is regular
if and only if it is isomorphic to a direct sum of band modules.
\end{Cor}

Let $B_1, \cdots, B_m$ be bands.
For positive integers $p_1, \cdots, p_m$ set
\begin{align*}
p &= \sum_{i=1}^m p_i,\\
n &= \sum_{i=1}^m p_i|B_i|,\\
F^p &= \{ (\lambda_1, \cdots, \lambda_p) \in \field^p 
\mid \lambda_i \not= \lambda_j \not= 0
\text{ for all } i \not= j \}.
\end{align*}
Define a morphism of varieties
\[ 
\GL_n(\field) \times F^p
\longrightarrow \V(n,a,b) 
\]
\[
(g,(\lambda_{11}, \cdots, \lambda_{p_11}, \cdots,
\lambda_{1m}, \cdots, \lambda_{p_mm})) \mapsto
g \cdot \left( \bigoplus_{j=1}^m 
M(B_j,\lambda_{1j}, \cdots, \lambda_{p_jj}) \right).
\] 
The image of this morphism is denoted by 
\[
{\mathcal F} = {\mathcal F}((B_1,p_1),\cdots,(B_m,p_m)).
\] 
We say that ${\mathcal F}$ is a $p$-{\it parametric family}.
In case $p_i = 1$ for some $i$, we write also just $B_i$ instead of 
$(B_i,p_i)$.
It follows from \cite{Kr} that $\dm \orb(y)$ is constant for 
all $y$ in a given family ${\mathcal F}$.
The following lemma is straightforward.

\begin{Lem}\label{dimensionoffamilies}
Any $p$-parametric family
${\mathcal F}$ is constructible, irreducible and has dimension
$p + \dm \orb(y)$ where
$y$ is any point in ${\mathcal F}$.
\end{Lem}

\begin{Lem}
Each direct sum of band modules is contained in the closure of
some family ${\mathcal F}$.
\end{Lem}

\begin{proof}
A band module $M(B,\lambda,n)$ is obviously contained in the closure of
the set of all band modules $M(B,\lambda_1, \cdots, \lambda_n)$ where
the $\lambda_i$ are pairwise different.
\end{proof}

For a string $C$ define 
\[
{\mathcal P}(C) = \{ (D,E,F) \mid D,E,F \in \sub \text{ and } DEF = C \}.
\]
We call $(D,E,F) \in {\mathcal P}(C)$ a {\it factor string} of $C$ if the 
following hold:
\begin{enumerate}

\item[{\rm (1)}] Either $D=1$ or $D = d_1 \cdots d_n$ where
$d_n = x$;

\item[{\rm (2)}] Either $F=1$ or $F = f_1 \cdots f_m$ where
$f_1 = y$.

\end{enumerate}
Dually, we call $(D,E,F)$ a {\it substring} of $C$ if the 
following hold:
\begin{enumerate}

\item[{\rm (1)}] Either $D=1$ or $D = d_1 \cdots d_n$ where
$d_n = y$;

\item[{\rm (2)}] Either $F=1$ or $F = f_1 \cdots f_m$ where
$f_1 = x$.

\end{enumerate}
Let ${\rm fac}(C)$ be the set of factor strings of $C$, 
and by ${\rm sub}(C)$ we denote the set of substrings of $C$.
For strings $C_1$ and $C_2$ let
\[
{\mathcal A}(C_1,C_2) = \{ ((D_1,E_1,F_1),(D_2,E_2,F_2)) \in 
{\rm fac}(C_1) \times {\rm sub}(C_2) \mid E_1 = E_2 \}.
\]
For example, if $C_1 = xxy$ and $C_2 = xyxx$, then 
\begin{multline*}
{\mathcal A}(C_1,C_2) = \{ ((xx,1,y),(1,1,xyxx)),
((xx,1,y),(xy,1,xx)),\\
((1,xx,y),(xy,xx,1)),
((x,x,y),(xy,x,x)),
((x,xy,1),(1,xy,xx)) \}.
\end{multline*} 
For each $a = ((D_1,E_1,F_1),(D_2,E_2,F_2)) \in {\mathcal A}(C_1,C_2)$ 
we define a homomorphism 
\[ 
f_a : M(C_1) \longrightarrow M(C_2) 
\] 
as follows:
Define
\[
f_a(z_{|D_1|+i}) = z_{|D_2|+i}
\] 
for $1 \leq i \leq |E_1|+1$, 
and all other canonical basis vectors of $M(C_1)$ are mapped to 0. 
Such homomorphisms are called {\it graph maps}.
The following theorem is a special case of the main result in \cite{CB}.

\begin{Thm}[Crawley-Boevey]\label{graphmaps}
The graph maps $\{ f_a \mid a \in {\mathcal A}(C_1,C_2) \}$ form a 
$\field$-basis
of the homomorphism space $\Hom_\Lambda(M(C_1),M(C_2))$.
\end{Thm}

There is the following multiplicative behaviour of graph maps:
Let $f_a: M(C_1) \to M(C_2)$ and 
$f_b: M(C_2) \to M(C_3)$ be graph maps.
Then the composition $f_a f_b: M(C_1) \to M(C_3)$ is either 0
or a graph map.


\section{Index modules of Gelfand-Ponomarev algebras}\label{indexsection}


A module $M$ is called {\it biserial} if it is isomorphic to
\[
\bigoplus_{i=1}^m M(x^iy^j)
\]
where $0 \leq i \leq a-1$ and $0 \leq j \leq b-1$.
For example, $\Lambda$ regarded as a module over itself is isomorphic 
to the biserial module $M(x^{a-1}y^{b-1})$.
Note also that any projective $\Lambda$-module is isomorphic to
$\Lambda^n$ for some $n \geq 1$.

\begin{Lem}\label{stringalgebrasaresubfinite}
Gelfand-Ponomarev algebras are subfinite, and all their
index modules are biserial.
\end{Lem}

\begin{proof}
Any 
Gelfand-Ponomarev algebra $\Lambda$ is a monomial algebra.
Thus by \cite[Lemma 3]{ZH}, if $U$ is a submodule of a projective 
$\Lambda$-module, then
\[
(U \cdot x) \cap (U \cdot y) = 0.
\]
It follows from the description of indecomposable $\Lambda$-modules 
that the biserial modules are the only $\Lambda$-modules 
which have this property.
\end{proof}

A case by case analysis shows the following:

\begin{Lem}\label{indexmodules}
A biserial $\Lambda$-module  
\[ 
L = S^{m_s} \oplus \bigoplus\limits_{i \geq 1} M(x^i)^{m_{xi}} \oplus 
\bigoplus\limits_{j \geq 1} M(y^j)^{m_{yj}} 
\oplus \bigoplus_{i,j \geq 1} M(x^iy^j)^{m_{ij}} 
\]
is isomorphic to a submodule of $\Lambda^n$ if and only if 
the following hold:
\begin{eqnarray*}
m_{xi} \not= 0   & \Longrightarrow & i \leq a-2,\\
m_{yj} \not= 0   & \Longrightarrow & j \leq b-2,\\
m_{ib-1} \not= 0 & \Longrightarrow & i = a-1,\\
m_{a-1j} \not= 0 & \Longrightarrow & j = b-1,\\
\sum_{i \geq 1} m_{xi} + \sum_{i,j \geq 1} m_{ij} &\leq& n,\\
\sum_{j \geq 1} m_{yj} + \sum_{i,j \geq 1} m_{ij} &\leq& n,\\
m_s+ \sum_{i \geq 1} m_{xi} + \sum_{j \geq 1} m_{yj} + 
2 \left( \sum_{i,j \geq 1} m_{ij} \right) &\leq& 2n.
\end{eqnarray*}
The dimension of $L$ is 
\[ 
m_s + \sum_{i \geq 1} m_{xi}(i+1) + 
\sum_{j \geq 1} m_{yj}(j+1) + 
\sum_{i,j \geq 1} m_{ij}(i+j+1). 
\]
\end{Lem}

\begin{Lem}\label{dimhomindextoproj}
Let $L \in \imod_\Lambda(n)$ and assume that $L$ is the direct sum of $m$
indecomposable modules.
Let $p$ be the number of indecomposable 
projective modules among these direct summands.
Then we have
\[ 
\dm \Hom(L,\Lambda) = n(d-1) + m - p 
\]
where $d = {\rm dim}(\Lambda)$.
\end{Lem}

\begin{proof}
We have $\dm \Hom(\Lambda,\Lambda) = \dim(\Lambda) = d$ and so 
$\dm \Hom(P,\Lambda) = \dm(P)$ for any projective module $P$.
Each indecomposable non-projective direct summand of $L$ is of the form
$M(x^iy^j)$ with $0 \leq i \leq a-2$ and $0 \leq j \leq b-2$.
We have
\begin{align*}
\dm M(x^iy^j) &= i+j+1,\\
\dm \Hom(M(x^iy^j),\Lambda) &= i+j+2. 
\end{align*}
This can be checked directly or by applying \cite{CB}.
Since ${\rm dim}(L) = n(d-1)$, the result follows.
\end{proof}

\begin{Lem}\label{flipdegeneration}
If $0 \leq p \leq i \leq a-1$ and 
$0 \leq q \leq j \leq b-1$, then
\[
M(x^iy^q) \oplus M(x^py^j) \leq_{\rm deg} M(x^iy^j) \oplus
M(x^py^q).
\]
\end{Lem}

\begin{proof}
One can easily construct a short exact sequence
\[ 
0 \longrightarrow M(x^iy^j) \longrightarrow M(x^iy^q) \oplus 
M(x^py^j) \longrightarrow M(x^py^q) \longrightarrow 0. 
\]
\end{proof}

A degeneration of the same form as in the previous lemma is called
a {\it flip degeneration}. (We `flip' $q$ and $j$.)
An index module $L$ is called {\it flip minimal} if it is isomorphic
to a direct sum of the form
\[ 
\Lambda^p \oplus \bigoplus\limits_{i=1}^t M(x^{c_i}y^{d_{t-i+1}}) 
\]
such that $c_i \geq c_{i+1}$, $d_{i} \geq d_{i+1}$, $0 \leq c_i \leq a-2$ 
and $0 \leq d_i \leq b-2$ for all $i$.
It follows from the previous two lemmas that for any
index module $L$ there exists a chain 
\[ 
L_1 \leq_{\rm deg} L_2 \leq_{\rm deg} \cdots \leq_{\rm deg} L_t = L 
\]
of flip degenerations of index modules with $L_1$ being flip minimal and 
\[ 
\dm \Hom(L_i,\Lambda) =  \dm \Hom(L_1,\Lambda) 
\]
for all $i$.

\begin{Lem}\label{boxmovedegeneration}
If $1 \leq p \leq i \leq a-2$ and 
$0 \leq q,j \leq b-1$, then
\[
M(x^{i+1}y^j) \oplus M(x^{p-1}y^q) \leq_{\rm deg} M(x^iy^j)
\oplus M(x^py^q).
\] 
\end{Lem}

\begin{proof}
One can construct a short exact sequence
\[ 
0 \longrightarrow M(x^py^q) \longrightarrow M(x^{i+1}y^j) \oplus 
M(x^{p-1}y^q) \longrightarrow M(x^iy^j) \longrightarrow 0. 
\]
\end{proof}

A degeneration of the same form as in the above lemma is called a 
{\it box move degeneration}. 
The modules over $\field[x]/(x^n)$ correspond to partitions, 
or equivalently to Young diagrams, and the degenerations of these modules
are given by moving boxes of the Young diagrams.
We are in a similar situation here.
Note that Lemma \ref{boxmovedegeneration} has an obvious dual version,
exchanging the roles of $x$ and $y$.


\section{Regular strata are irreducible}\label{sectionregcomp}

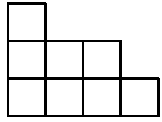
\begin{figure}[ht]
\unitlength0.5cm
\begin{picture}(4,3)
\put(0,0){\line(0,1){3.0}}
\put(1,0){\line(0,1){3.0}}
\put(2,0){\line(0,1){2.0}}
\put(3,0){\line(0,1){2.0}}
\put(4,0){\line(0,1){1.0}}
\put(0,0){\line(1,0){4.0}}
\put(0,1){\line(1,0){4.0}}
\put(0,2){\line(1,0){3.0}}
\put(0,3){\line(1,0){1.0}}
\end{picture}
\caption{The Young diagram $Y(3,2,2,1)$}\label{fig3}
\end{figure}

For a partition ${\bf p} = (p_1, \cdots, p_t)$ let
$Y({\bf p})$ be its corresponding Young diagram, which
has $p_i$ boxes in the $i$th column.
For example, the Young diagram
$Y(3,2,2,1)$ looks as in Figure \ref{fig3}.
For a partition ${\bf p} = (p_1, \cdots,p_t)$ the {\it dual partition}
is defined as
${\bf p}^* = (r_1, \cdots, r_{p_1})$, 
where the $r_j$ are the number of boxes in the rows of the
Young diagram $Y({\bf p})$, ordered decreasingly.
For example, 
\[
(3,2,2,1)^* = (4,3,1).
\]
Now let $A \in \nil(n,a)$ with ${\bf p} = p(A)$.
Then the boxes of the Young diagram $Y({\bf p})$ can be considered
as a certain basis of $\field^n$, and $A$ can be considered as an 
endomorphism of $\field^n$.
If $b$ is a box such that there is a box $b'$ below $b$, then $A$
maps $b$ to $b'$, and $b$ is mapped to 0, otherwise.
Figure \ref{fig3.1} illustrates this for $p(A) = (3,2,2,1)$, 
where the arrows indicate how $A$ acts on the boxes.
\begin{figure}[ht]
\unitlength0.5cm
\begin{picture}(4,4)
\put(0,1.2){\line(0,1){3.0}}
\put(1,1.2){\line(0,1){3.0}}
\put(2,1.2){\line(0,1){2.0}}
\put(3,1.2){\line(0,1){2.0}}
\put(4,1.2){\line(0,1){1.0}}
\put(0,1.2){\line(1,0){4.0}}
\put(0,2.2){\line(1,0){4.0}}
\put(0,3.2){\line(1,0){3.0}}
\put(0,4.2){\line(1,0){1.0}}
\put(0.5,3.6){\vector(0,-1){0.9}}
\put(0.5,2.6){\vector(0,-1){0.9}}
\put(1.5,2.6){\vector(0,-1){0.9}}
\put(2.5,2.6){\vector(0,-1){0.9}}
\put(0.5,1.6){\vector(0,-1){0.9}}
\put(1.5,1.6){\vector(0,-1){0.9}}
\put(2.5,1.6){\vector(0,-1){0.9}}
\put(3.5,1.6){\vector(0,-1){0.9}}
\put(0.3,0){$0$}
\put(1.3,0){$0$}
\put(2.3,0){$0$}
\put(3.3,0){$0$}
\end{picture}
\caption{}\label{fig3.1}
\end{figure}
Now let ${\bf p}^* = (r_1, \cdots, r_m)$ be the dual partition
of ${\bf p}$.
Then $r_1 = \dm \ke(A)$, $r_1 + r_2 = \dm \ke(A^2)$ etc.

If ${\bf p}$ and ${\bf q}$ are arbitrary partitions, then define
${\bf p} \trianglelefteq {\bf q}$ if 
\[ 
\sum_{i=1}^l p_i \leq \sum_{j=1}^l q_j 
\]
for all $l$, where we set $p_i = 0$ and $q_j = 0$ for
all $i > l({\bf p})$ and $j > l({\bf q})$.
This partial order is usually called the {\it dominance order}. 
The proof of the following proposition can be found 
in \cite{Ge}, see also \cite{H}.

\begin{Prop}\label{nilvariety}
For ${\bf p} \in {\mathcal P}(n,l)$ we have
\[ 
\dm C({\bf p}) = n^2- \sum_{i=1}^t r_i^2 
\]
where ${\bf p}^* = (r_1, \cdots, r_t)$,
\[ 
C({\bf p}) = \{ A \in \nil(n,l) \mid \rk(A^k) = n- \sum_{j=1}^k r_j,
1 \leq k \leq t \} 
\] 
and 
\[ 
\overline{C({\bf p})} =
\{ A \in \nil(n,l) \mid \rk(A^k) \leq n- \sum_{j=1}^k r_j,
1 \leq k \leq t \}. 
\] 
In particular,
if ${\bf p},{\bf q} \in {\mathcal P}(n,l)$,
then $C({\bf p}) \subseteq \overline{C({\bf q})}$
if and only if ${\bf p} \trianglelefteq {\bf q}$.
\end{Prop}

Recall that we defined two maps
\[
\xymatrix{
 & \V(n,a,b) \ar[dr]^{\pi_2} \ar[dl]_{\pi_1} \\
\nil(n,a) & & \nil(n,b)
}
\]
with $\pi_1(A,B) = A$ and $\pi_2(A,B) = B$, 
and for $({\bf a},{\bf b}) \in {\mathcal P}(n,a,b)$ we set
\begin{align*}
\Delta({\bf a}) &= \pi_1^{-1}(C({\bf a})),\\
\Delta({\bf a},{\bf b}) &= \pi_1^{-1}(C({\bf a})) \cap 
\pi_2^{-1}(C({\bf b})).
\end{align*}
Thus, as a consequence of Proposition \ref{nilvariety} we get 
\[
\Delta({\bf a}) = 
\{ (A,B) \in \V(n,a,b) \mid 
\rk(A^k) = n-\sum_{j=1}^k m_j, 1 \leq k \leq r \}
\]
and
\begin{multline*} 
\Delta({\bf a},{\bf b}) = 
\{ (A,B) \in \V(n,a,b) \mid 
\rk(A^k) = n-\sum_{j=1}^k m_j, 1 \leq k \leq r,\\ 
\rk(B^l) = n-\sum_{j=1}^l n_j,
1 \leq l \leq s \}
\end{multline*}
where ${\bf a}^* = (m_1, \cdots, m_r)$ and
${\bf b}^* = (n_1, \cdots, n_s)$.
In particular, $\Delta({\bf a})$ and  
$\Delta({\bf a},{\bf b})$ are locally closed in $\V(n,a,b)$.

The following lemma is an easy exercise.

\begin{Lem}\label{basic1}
If $M \in \nil(n,l)$, then
$\rk(M) = n - l(p(M))$.
\end{Lem}

\begin{proof}[Proof of Proposition \ref{regularstrata}]
Let $(A,B) \in \V(n,a,b)$, and set
\[
({\bf a},{\bf b}) = (p(A),p(B)) \in {\mathcal P}(n,a,b).
\]
Assume that
$(A,B)$ is regular, i.e. $\rk(A) + \rk(B) = n$.
By Lemma \ref{basic1} this is equivalent to $n = l({\bf a}) + l({\bf b})$.
Thus, if $\Delta({\bf a},{\bf b})$ contains a regular element, then all
elements in $\Delta({\bf a},{\bf b})$ are regular.
We know that $(A,B)$
is isomorphic to a direct sum of band modules.
But any band is (up to equivalence) of the form 
$x^{c_1}y^{d_1} \cdots x^{c_t}y^{d_t}$ with $c_i,d_i \geq 1$ for
$1 \leq i \leq t$.
This implies that the number of entries which are at least 2 in ${\bf a}$
is equal to the number of entries which are at least 2 in ${\bf b}$.
In other words, $l({\bf a}-1) = l({\bf b}-1)$.
Conversely, if $({\bf a},{\bf b}) \in {\mathcal P}(n,a,b)$ with
$l({\bf a}) + l({\bf b}) = n$, $l({\bf a}-1) = l({\bf b}-1)$, 
${\bf a}-1 = (c_1, \cdots, c_t)$ and ${\bf b}-1 = (d_1, \cdots, d_t)$,
then set 
\[
(A,B) = M(x^{c_1}y^{d_1} \cdots x^{c_t}y^{d_t},\lambda).
\]
Clearly, we have $p(A) = {\bf a}$, $p(B) = {\bf b}$,
and $(A,B)$ (and therefore also $({\bf a},{\bf b})$) is regular.
This finishes the proof.
\end{proof}

Altogether, we get that
for a regular $(A,B) \in \V(n,a,b)$ the following are equivalent:
\begin{itemize}

\item $\dm \tp(A,B) = p$;

\item $\dm \soc(A,B) = p$;

\item $\dm (\ke(A) \cap \ke(B)) = p$;

\item $\dm (\im(A) \cap \im(B)) = p$;

\item $l(p(A)-1) = p$;

\item $l(p(B)-1) = p$.

\end{itemize}
{\bf Example}:
Let $(A,B) = M(xxyxy,\lambda)$.
Then $p(A) = (3,2)$ and $p(B) = (2,2,1)$.
Thus, $l(p(A)-1) = l(2,1) = 2$ and
$l(p(B)-1) = l(1,1) = 2$.
It is also clear that $M(xxyxy,\lambda)$ has a 2-dimensional socle and
a 2-dimensional top.
As an illustration, see Figure \ref{fig4}.
\begin{figure}[ht]
\unitlength1.0cm
\begin{picture}(4,3)

\put(0.5,0.5){$z_{11}$}
\put(3.3,0.6){\vector(-1,0){2.2}}\put(1.8,0.3){$(y,\lambda)$}
\put(3.7,0.8){\vector(0,1){1.0}}\put(3.8,1.2){$x$}
\put(3.5,0.5){$z_{51}$}

\put(0.5,2){$z_{21}$}
\put(0.7,1.8){\vector(0,-1){1.0}}\put(0.4,1.2){$x$}
\put(1.8,2.1){\vector(-1,0){0.7}}\put(1.5,2.2){$x$}
\put(2,2){$z_{31}$}
\put(2.6,2.1){\vector(1,0){0.7}}\put(2.8,2.2){$y$}
\put(3.5,2){$z_{41}$}

\end{picture}
\caption{The band module $M(xxyxy,\lambda)$}\label{fig4}
\end{figure}

The next lemma follows directly from the construction of projective
covers of indecomposable $\Lambda$-modules.
These covers are easy to construct.

\begin{Lem}\label{indexsummands}
Assume that $\sub(L) \subset \V(n,a,b)$ contains a regular element $(A,B)$.
Then $L$ is a direct sum of $n$ indecomposable modules, and exactly 
$n - \dm \tp(A,B)$ of these are projective.
\end{Lem}

A $\Lambda$-module is called a  {\it diamond module} if it is isomorphic
to $M(x^iy^j,\lambda)$ for some $1 \leq i \leq a-1$ and
$1 \leq j \leq b-1$.
Thus the diamond modules are the band modules with simple top (and
therefore also with simple socle).
We now associate to any
regular element $({\bf a},{\bf b})$ a {\it diamond family} 
${\mathcal F}({\bf a},{\bf b})$
which consists of direct sums of diamond modules. 

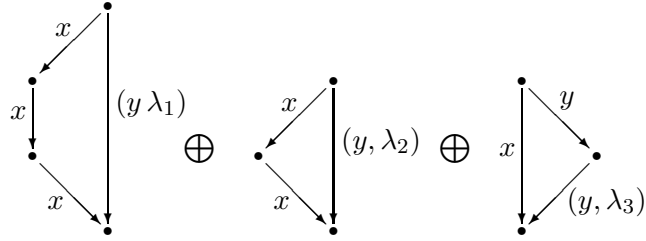
\begin{figure}[ht]
\unitlength1.0cm
\begin{picture}(9,4)

\put(2,0.5){\circle*{0.1}}
\put(1,1.5){\circle*{0.1}}
\put(1,2.5){\circle*{0.1}}
\put(2,3.5){\circle*{0.1}}
\put(1.9,3.4){\vector(-1,-1){0.8}}
\put(1.3,3.1){$x$}
\put(1,2.4){\vector(0,-1){0.8}}
\put(0.7,2.0){$x$}
\put(1.1,1.4){\vector(1,-1){0.8}}
\put(1.2,0.8){$x$}
\put(2,3.4){\vector(0,-1){2.8}}
\put(2.1,2.1){$(y\,\lambda_1)$}

\put(3,1.5){$\bigoplus$}

\put(5,0.5){\circle*{0.1}}
\put(4,1.5){\circle*{0.1}}
\put(5,2.5){\circle*{0.1}}
\put(4.9,2.4){\vector(-1,-1){0.8}}
\put(4.3,2.1){$x$}
\put(4.1,1.4){\vector(1,-1){0.8}}
\put(4.2,0.8){$x$}
\put(5,2.4){\vector(0,-1){1.8}}
\put(5.1,1.6){$(y,\lambda_2)$}

\put(6.4,1.5){$\bigoplus$}

\put(7.5,0.5){\circle*{0.1}}
\put(8.5,1.5){\circle*{0.1}}
\put(7.5,2.5){\circle*{0.1}}
\put(7.5,2.4){\vector(0,-1){1.8}}
\put(7.2,1.5){$x$}
\put(7.6,2.4){\vector(1,-1){0.8}}
\put(8.0,2.2){$y$}
\put(8.4,1.4){\vector(-1,-1){0.8}}
\put(8.1,0.8){$(y,\lambda_3)$}

\end{picture}
\caption{An element in ${\mathcal F}((4,3,2,1),(3,2,2,1,1,1))$}\label{fig5}
\end{figure}

Let $({\bf a},{\bf b}) \in {\mathcal P}(n,a,b)$ be regular.
Thus, $l({\bf a}-1) = l({\bf b}-1)$ by Proposition \ref{regularstrata}.
Assume that ${\bf a}-1 = (c_1, \cdots, c_t)$ and 
${\bf b}-1 = (d_1, \cdots, d_t)$.
Let 
\[ 
{\mathcal F}({\bf a},{\bf b}) =
{\mathcal F}(x^{c_1}y^{d_t},x^{c_2}y^{d_{t-1}}, \cdots, x^{c_t}y^{d_1}). 
\]
Thus every module in ${\mathcal F}({\bf a},{\bf b})$ is isomorphic to
\[
\bigoplus_{i=1}^t M(x^{c_i}y^{d_{t-i+1}},\lambda_i)
\]
for some pairwise different $\lambda_i$.
For example, a module in 
\[
{\mathcal F}((4,3,2,1),(3,2,2,1,1,1))
\]
looks as in Figure \ref{fig5}, where the points 
are just the basis vectors of the module.
Note that ${\mathcal F}({\bf a},{\bf b}) \subset \Delta({\bf a},{\bf b})$.

\begin{Prop}\label{reductiontodiamonds}
If $({\bf a},{\bf b}) \in {\mathcal P}(n,a,b)$ is regular, then
${\mathcal F}({\bf a},{\bf b})$ is dense in 
$\Delta({\bf a},{\bf b})$ and has dimension
\[ 
n^2 - \sum_{i=1}^r m_i^2 - \sum_{i=1}^s n_i^2 +l({\bf a}-1)^2 
\]
where $({\bf a}-1)^* = (m_1, \cdots, m_r)$ and 
$({\bf b}-1)^* = (n_1, \cdots, n_s)$.
In particular, $\Delta({\bf a},{\bf b})$ is irreducible.
\end{Prop}

\begin{proof}
Let $(A,B) \in \Delta({\bf a},{\bf b})$ be regular.
Thus $(A,B)$ is in some stratum $\sub(L)$ with $L$
a direct sum of $n$ indecomposable modules, and exactly 
$n-l({\bf a}-1)$ of them are projective, see Lemma \ref{indexsummands}.
By Lemma \ref{dimhomindextoproj} we get 
\[ 
\dm \Hom(L,\Lambda) = n(d-1) + l({\bf a}-1) 
\]
where $d = {\rm dim}(\Lambda)$.
Assume ${\bf a}-1 = (c_1, \cdots, c_t)$ and ${\bf b}-1 = (d_1, \cdots, d_t)$.
By Proposition \ref{regularstrata}
each module in $\Delta({\bf a},{\bf b})$ is isomorphic to a 
direct sum of band modules, and
one checks easily that $L = (E,F)$ with 
\[
p(E) = (e_1, \cdots, e_t) = (a-c_t-1,a-c_{t-1}-1, \cdots, a-c_1-1)
\]
and
\[
p(F) = (f_1, \cdots, f_t) = (b-d_t-1,b-d_{t-1}-1, \cdots, b-d_1-1).
\] 
Define
\[ 
L({\bf a},{\bf b}) = \Lambda^{n-t} \oplus 
\bigoplus\limits_{i=1}^t M(x^{e_i}y^{f_{t-i+1}}). 
\]
We apply a sequence of flip 
degenerations to $L$ and get
\[
L({\bf a},{\bf b}) \leq_{\rm deg} L
\]
with
\[
\dm \Hom(L({\bf a},{\bf b}),\Lambda) = \dm \Hom(L,\Lambda).
\]
Then Theorem \ref{richmond},(3) yields
that the stratum $\sub(L({\bf a},{\bf b}))$ is dense in 
$\Delta({\bf a},{\bf b})$.
By \ref{richmond},(1) we get that $\Delta({\bf a},{\bf b})$ is
irreducible.
Observe that 
\[
{\mathcal F}({\bf a},{\bf b}) \subset
\sub(L({\bf a},{\bf b})).
\]
We have
\[
\dm {\mathcal F}({\bf a},{\bf b}) =
n^2 - \sum_{i=1}^r m_i^2 - \sum_{i=1}^s n_i^2 +l({\bf a}-1)^2 
\]
where $({\bf a}-1)^* = (m_1, \cdots, m_r)$ and 
$({\bf b}-1)^* = (n_1, \cdots, n_s)$.
This follows from Lemma \ref{dimensionoffamilies}, the dimension formula
for orbits and \cite{Kr}.
Using the dimension formula in Theorem \ref{richmond},(1)
and applying \cite{CB} we get 
\[ 
\dm \sub(L({\bf a},{\bf b})) = \dm {\mathcal F}({\bf a},{\bf b}).
\]
This implies that ${\mathcal F}({\bf a},{\bf b})$ is dense in 
$\Delta({\bf a},{\bf b})$.
\end{proof}

Thus, from the above proposition we get the remarkable result that the 
diamond families form a dense subset in the set of all regular elements
in $\V(n,a,b)$.


\section{The nilpotent case}\label{nilcase1}


The following is easy to prove.

\begin{Lem}\label{nonreg3}
For $u \in \{ x,y \}$ and strings $C$ and $D$
the following hold:
\begin{itemize}
\item[(1)]
If $CuD$ is a string, then
\[
M(C) \oplus M(D) \in \overline{\orb(M(CuD))};
\]

\item[(2)]
If $Cu$ is a band, then
\[
M(C) \in \overline{{\mathcal F}(Cu)} =  \overline{{\mathcal F}(uC)}.
\]

\end{itemize}
\end{Lem}

\begin{Lem}\label{stringreduction}
If $1 \leq i \leq a-1$, $1 \leq j \leq b-1$, and $l \geq 0$ such
that $j+l+1 \leq b-1$, then
\[
M(x^iy^j,\lambda) \oplus M(y^l) \in
\overline{{\mathcal F}(x^iy^{j+l+1})}.
\]
\end{Lem}

\begin{proof}
There exists a short exact sequence
\[ 
0 \longrightarrow M(x^iy^j,\lambda) \longrightarrow M(y^{j+l}x^i) 
\longrightarrow M(y^l) \longrightarrow 0. 
\]
Thus
\[
M(y^{j+l}x^i) \leq_{\rm deg} M(x^iy^j,\lambda) \oplus M(y^l).
\]
Then we use Lemma \ref{nonreg3},(2).
\end{proof}

\begin{Lem}\label{stringreduction3}
Let $(C,D) \in \V(n,n,n)$ with $\rk(C) + \rk(D) < n$ and
$\rk(D) < n-1$.
Then $(C,D)$ is contained in the closure of
\[
\{ (A,B) \in \V(n,n,n) \mid \rk(A) = \rk(C), \rk(B) = \rk(D)+1 \}.
\]
\end{Lem}

\begin{proof}
Set
\[
{\mathcal C} =
\{ (A,B) \in \V(n,n,n) \mid \rk(A) = \rk(C), \rk(B) = \rk(D)+1 \}, 
\]
and let $s = n-\rk(C)-\rk(D)$.
Thus $(C,D)$ is isomorphic to a module
\[
M \oplus \bigoplus_{i=1}^s M(C_i)
\]
where $M = 0$ or $M$ is a direct sum of band modules.
There are three cases to consider:
First, if $s \geq 2$, then 
\[
(C,D) \in \overline{\orb(M \oplus M(C_1yC_2) \oplus 
M(C_3) \oplus \cdots \oplus M(C_s))} 
\subseteq \overline{{\mathcal C}}.
\]
Second, if $s=1$ and $C_1 \not= y^l$ for some $l \geq 0$, then
\[
(C,D) \in \overline{M \oplus {\mathcal F}(C_1y)}
\subseteq \overline{{\mathcal C}}.
\]
Finally, assume that $s=1$ and $C_1 = y^l$ for some $l \geq 0$.
Since $\rk(D) < n-1$, this implies $l < n-1$ and thus $M \not= 0$.
Using Proposition \ref{reductiontodiamonds} 
we can assume without loss of generality that
$M$ is a direct sum of diamond modules. 
Let $M(x^iy^j,\lambda)$ be one of these direct summands,
thus $M = M' \oplus M(x^iy^j,\lambda)$ for some $M'$.
Then we use Lemma \ref{stringreduction} and get
\[
(C,D) \in \overline{M' \oplus {\mathcal F}(x^iy^{j+l+1})}
\subseteq \overline{{\mathcal C}}.
\]
Note that we used several times our assumption $a=b=n$ by assuming
that certain words in $x$ and $y$ are actually strings, i.e. that they
do not contain subwords of the form $x^a$ or $y^b$.
This finishes the proof.
\end{proof}

\begin{Cor}\label{stringreduction2}
Let $(A,B) \in \V(n,n,n)$ with 
$\rk(A) \leq n-i$ and $\rk(B) \leq i$.
Then $(A,B)$ is contained in the closure of
\[
\{ (A,B) \in \V(n,n,n) \mid \rk(A) = n-i, \rk(B) = i \}.
\]
\end{Cor}

\begin{Lem}\label{degoftopsocleonemodules}
If $u_1,u_2,v_1,v_2 \geq 1$, $u_1+u_2 \leq a-1$ and $v_1+v_2 \leq b-1$, then 
\[ 
M(x^{u_1+u_2}y^{v_1+v_2},-\lambda_1\lambda_2) \leq_{\rm deg}
M(x^{u_1}y^{v_1},\lambda_1) \oplus M(x^{u_2}y^{v_2},\lambda_2). 
\]  
\end{Lem}

\begin{proof}
It is straightforward to construct a short exact sequence
\[ 
0 \longrightarrow  M(x^{u_1}y^{v_1},\lambda_1) \longrightarrow
M(x^{u_1+u_2}y^{v_1+v_2}, -\lambda_1\lambda_2) \longrightarrow 
M(x^{u_2}y^{v_2},\lambda_2) \longrightarrow 0. 
\]
\end{proof}

\begin{proof}[Proof of Theorem \ref{nil2}]
Let $({\bf a},{\bf b}) \in {\mathcal P}(n,n,n)$ be regular with
\[
{\bf a}-1 = (c_1, \cdots, c_t)
\]
and
\[
{\bf b}-1 = (d_1, \cdots, d_t).
\]
The diamond family ${\mathcal F}({\bf a},{\bf b})$ is dense
in $\Delta({\bf a},{\bf b})$ by Proposition \ref{reductiontodiamonds},
and each module in ${\mathcal F}({\bf a},{\bf b})$ is of the form
\[
\bigoplus_{j=1}^t M(x^{c_j}y^{d_{t-j+1}},\lambda_j)
\]
for some $\lambda_j$.
Since $a = b = n$, we know that $x^{n-i}y^i$ is a string for all
$1 \leq i \leq n-1$.
Now we use Lemma \ref{degoftopsocleonemodules} and get that
\[
{\mathcal F}({\bf a},{\bf b}) \subset \Delta({\bf a},{\bf b})
\subset \overline{{\mathcal F}(x^{n-i}y^i)},
\]
where 
\[
n-i = \sum_{j=1}^t c_j
\]
and
\[
i = \sum_{j=1}^t d_j.
\]
This implies 
\[
\{ (A,B) \in \V(n,n,n) \mid \rk(A) = n-i, \rk(B) = i \} \subset 
\overline{{\mathcal F}(x^{n-i}y^i)}. 
\]
Then Corollary \ref{stringreduction2} implies 
\[
\overline{{\mathcal F}(x^{n-i}y^i)} =  
\{ (A,B) \in \V(n,n,n) \mid \rk(A) \leq n-i, \rk(B) \leq i \}.
\]
By Proposition \ref{reductiontodiamonds} we get
\[
\dm {\mathcal F}(x^{n-i}y^i) = \dm \overline{{\mathcal F}(x^{n-i}y^i)} =   
n^2-n+1.
\]
This finishes the proof.
\end{proof}


\section{Classification of regular irreducible components}
\label{classregcomp}


If $({\bf a},{\bf b}), ({\bf c},{\bf d}) \in {\mathcal P}(n,a,b)$ 
with ${\bf a} \trianglelefteq {\bf c}$ and ${\bf b} \trianglelefteq 
{\bf d}$, then we
write $({\bf a},{\bf b}) \trianglelefteq ({\bf c},{\bf d})$.
This defines a partial order on ${\mathcal P}(n,a,b)$.

\begin{Lem}\label{regpairs}
If $({\bf a},{\bf b}), ({\bf c},{\bf d}) \in {\mathcal P}(n,a,b)$ 
are regular with
$({\bf a},{\bf b}) \trianglelefteq ({\bf c},{\bf d})$, then
$l({\bf a}) = l({\bf c})$ and $l({\bf b}) = l({\bf d})$.
\end{Lem}

\begin{proof}
For all regular pairs $({\bf e},{\bf f})$ we
have $l({\bf e}) + l({\bf f}) = n$.
Since ${\bf a} \trianglelefteq {\bf c}$, we have $l({\bf a}) \geq 
l({\bf c})$,
and from ${\bf b} \trianglelefteq {\bf d}$ we get $l({\bf b}) \geq 
l({\bf d})$.
This implies $l({\bf a}) = l({\bf c})$ and $l({\bf b}) = l({\bf d})$.
\end{proof}

The next lemma is a consequence of Proposition \ref{nilvariety}.

\begin{Lem}\label{strataintersection}
Let $({\bf a},{\bf b}), ({\bf c},{\bf d}) \in
{\mathcal P}(n,a,b)$.
If 
\[
\Delta({\bf a},{\bf b}) \cap \overline{\Delta({\bf c},{\bf d})} 
\not=  \emptyset, 
\] 
then $({\bf a},{\bf b}) \trianglelefteq ({\bf c},{\bf d})$.
\end{Lem}

Let 
\begin{multline*}
{\mathcal P}_{i,\reg}^p  = {\mathcal P}_{i,\reg}^p(n,a,b) = 
\{ ({\bf a},{\bf b}) \in
{\mathcal P}(n,a,b) \mid l({\bf a}) = i, l({\bf b}) = n-i,\\ 
l({\bf a}-1) = p \},
\end{multline*}
and
\[ 
\V_{i,\reg}^p = \V_{i,\reg}^p(n,a,b) = 
\bigcup_{({\bf a},{\bf b}) \in {\mathcal P}_{i,\reg}^p}
\Delta({\bf a},{\bf b}). 
\]
This implies
\begin{multline*}
\V_{i,\reg}^p = \{ (A,B) \in \V(n,a,b) \mid 
\rk(A) = n-i,\rk(B) = i,\\
\dm \tp(A,B) = p \}.
\end{multline*}
In particular, $\V_{i,\reg}^p(n,a,b)$ is locally closed.

\begin{Prop}\label{vpiisirreducible}
Let $({\bf a},{\bf b}), ({\bf c},{\bf d}) \in 
{\mathcal P}_{i,\reg}^p(n,a,b)$.
Then 
\[
\Delta({\bf a},{\bf b}) \subset \overline{\Delta({\bf c},{\bf d})}
\]
if and only if 
$({\bf a},{\bf b}) \trianglelefteq ({\bf c},{\bf d})$.
\end{Prop}

\begin{proof}
If $({\bf a},{\bf b}) \trianglelefteq ({\bf c},{\bf d})$ does not hold,
then we apply Lemma \ref{strataintersection} and get 
\[ 
\Delta({\bf a},{\bf b}) \cap \overline{\Delta({\bf c},{\bf d})} 
=  \emptyset. 
\] 
Next, assume that 
$({\bf a},{\bf b}) \trianglelefteq ({\bf c},{\bf d})$
holds.
By Lemma \ref{indexsummands} 
each element in $\V_{i,\reg}^p$ belongs to 
some stratum of the form $\sub(L)$ with $L$ a direct sum of $n$
indecomposables, and exactly $n - p$ of these are projective.
Since $({\bf a},{\bf b}) \trianglelefteq ({\bf c},{\bf d})$
and $({\bf a},{\bf b}), ({\bf c},{\bf d}) \in 
{\mathcal P}_{i,\reg}^p(n,a,b)$,
there exists a chain
\[ 
L({\bf c},{\bf d}) = L_1 \leq_{\rm deg} 
L_2 \leq_{\rm deg} \cdots \leq_{\rm deg} L_t = L({\bf a},{\bf c}) 
\]
of box move degenerations between index modules such that
$\dm \Hom(L_i,\Lambda)$ is constant for all $L_i$ in this chain.
Now we use the same arguments as in the proof of Proposition 
\ref{reductiontodiamonds}, and finally we apply 
Theorem \ref{richmond},(3).
This finishes the proof.
\end{proof}

An element 
$({\bf a},{\bf b}) \in {\mathcal P}_{i,\reg}^p(n,a,b)$
is called $(i,p)$-{\it maximal} if it is maximal in
${\mathcal P}_{i,\reg}^p(n,a,b)$ with respect to the partial
order $\trianglelefteq$.
Clearly, each non-empty ${\mathcal P}_{i,\reg}^p(n,a,b)$
contains a unique $(i,p)$-maximal element.

It follows easily that an element $({\bf a},{\bf b}) \in 
{\mathcal P}_{i,\reg}^p(n,a,b)$ is $(i,p)$-maximal if and only if 
the following hold:
\begin{itemize}

\item ${\bf a}$ has at most one entry different from $1$, $2$ and $a$;

\item ${\bf b}$ has at most one entry different from $1$, $2$ and $b$.

\end{itemize}
As a consequence of Propositions 
\ref{reductiontodiamonds} and \ref{vpiisirreducible} we get 
the following:

\begin{Cor}\label{cor7.4}
The set $\V_{i,\reg}^p(n,a,b)$ is locally closed and irreducible,
and if it is non-empty, then
it contains ${\mathcal F}({\bf a},{\bf b})$ as a dense subset, where
$({\bf a},{\bf b})$ is the unique $(i,p)$-maximal element in
${\mathcal P}_{i,\reg}^p(n,a,b)$.
\end{Cor}

\begin{proof}[Proof of Theorem \ref{theorem1}]
We characterize the $(i,p)$-maximal elements 
$({\bf a},{\bf b})$ such that the closure of $\Delta({\bf a},{\bf b})$
is an irreducible component.
By the preceding results, these are then all regular irreducible components.
Assume that $({\bf a},{\bf b})$ is $(i,p)$-maximal.
Thus
\[
{\bf a}-1 = ((a-1)^{p-r-1},a-v-1,1^r)
\] 
and
\[
{\bf b}-1 = ((b-1)^{p-s-1},b-w-1,1^s)
\] 
where 
$0 \leq v \leq a-2$, $0 \leq w \leq b-2$, $0 \leq r,s \leq p-1$,
$v=0 \Rightarrow r=0$ and $w=0 \Rightarrow s=0$.
By Corollary \ref{cor7.4} we have
\[
\overline{{\mathcal F}({\bf a},{\bf b})} = 
\overline{\Delta({\bf a},{\bf b})} =
\overline{\V_{i,\reg}^p}.
\]
We claim that 
the closure of ${\mathcal F}({\bf a},{\bf b})$ is an irreducible component
if and only if $r+s+1 \leq p$.

First, let $r+s+1 > p$.
This implies that there exist $u_1,u_2,v_1,v_2 \geq 1$ such that 
each module in 
${\mathcal F}({\bf a},{\bf b})$ has a direct summand isomorphic to
\[
M(x^{u_1}y^{v_1},\lambda_1) \oplus M(x^{u_2}y^{v_2},\lambda_2)
\] 
where $u_1+u_2 \leq a-1$ and $v_1+v_2 \leq b-1$.
Now we apply Lemma \ref{degoftopsocleonemodules} and see that 
${\mathcal F}({\bf a},{\bf b})$ is contained in the closure of some other 
family ${\mathcal F}({\bf c},{\bf d})$.
In particular, the closure of ${\mathcal F}({\bf a},{\bf b})$ cannot be an 
irreducible component.
This proves one direction of the statement.

Second, assume that $r+s+1 \leq p$.
Since the function $\rk(-)$ is lower semicontinuous, 
$\V_{i,\reg}^p$ cannot be contained in the
closure of some $\V_{j,\reg}^q$ with $i \not= j$.
It is also clear that  $\V_{i,\reg}^p$ cannot be in the
closure of $\V_{j,\reg}^q$ if $p < q $.
Because in that case, we have 
\[ 
\dm \Hom(M,S) = p < q = \dm \Hom(N,S) 
\] 
for all $M \in \V_{i,\reg}^p$ and all $N \in \V_{j,\reg}^q$.
This is a contradiction to the upper semicontinuity of
the function $\dm \Hom(-,S)$.

Thus, assume that $i = j$, $p > q$ and $r+s+1 \leq p$.
Then the dimension formula in Proposition \ref{reductiontodiamonds} 
yields 
\[
\dm \V_{i,\reg}^p \geq \dm V_{j,\reg}^q.
\]
Again this implies that $\V_{i,\reg}^p$ cannot be in the closure
of $\V_{j,\reg}^q$.
Thus the closure of $\V_{i,\reg}^p$ must be an irreducible component.
Finally, note that $l({\bf a}-1) \leq |a \in {\bf a}| + |b \in {\bf b}|
+1$ if and only if $r+s+1 \leq p$.
This finishes the proof.
\end{proof}

\begin{proof}[Proof of Theorem \ref{nil1}]
Let ${\bf a} \in {\mathcal P}(n,n)$ be a partition of $n$.
If ${\bf a} = (1, \cdots, 1)$, then
$\Delta({\bf a})$ is the union of the orbits of $n$-dimensional
modules of the form
\[
\bigoplus_{i \geq 0} M(y^i)^{m_i},
\] 
and $\orb(M(y^{n-1}))$ is dense in
$\Delta({\bf a})$.
Thus $\Delta({\bf a})$ is irreducible and 
\[
\Delta({\bf a}) \subset \overline{{\mathcal F}(xy^{n-1})}
= \overline{\Delta(2,1,\cdots,1)}. 
\]
Next, assume that ${\bf a} \not= (1, \cdots, 1)$.
Thus $l({\bf a}) = i$ for some $1 \leq i \leq n-1$.
Then there exists a unique maximal (with respect to $\trianglelefteq$) 
partition ${\bf a}^\circ$ such that $({\bf a},{\bf a}^\circ)$ is regular.
Namely, we have ${\bf a}^\circ = (r_1, \cdots, r_{n-i})$ where
\[
r_j =  
\begin{cases}
i-l({\bf a}-1)+2 & \text{if $j=1$},\\
2 & \text{if $2 \leq j \leq l({\bf a}-1)$},\\
1 & \text{otherwise}.
\end{cases}
\]
Here we use our assumption $a=b=n$.
By Proposition \ref{vpiisirreducible}, we know that 
for any regular element $({\bf a},{\bf d})$ we have
\[
\Delta({\bf a},{\bf d}) \subset \overline{\Delta({\bf a},{\bf a}^\circ)}.
\]
Now, assume that $({\bf a},{\bf c})$ is non-regular with
$\Delta({\bf a},{\bf c})$ non-empty.
It follows from Lemma \ref{stringreduction3} that
\[
\Delta({\bf a},{\bf c}) \subset \overline{\Delta({\bf a},{\bf d})}
\]
for some regular $({\bf a},{\bf d})$.

This proves that $\Delta({\bf a})$ has $\Delta({\bf a},{\bf a}^\circ)$
as a dense subset.
Thus $\Delta({\bf a})$ is irreducible.

Recall that for regular elements, 
$({\bf a},{\bf b}) \trianglelefteq ({\bf c},{\bf d})$ 
implies $l({\bf a}) = l({\bf c})$, see Lemma \ref{regpairs}.
Using Lemma \ref{strataintersection}, we get that
\[ 
\Delta({\bf a}) \subset \overline{\Delta({\bf c})} 
\] 
implies
${\bf a} \trianglelefteq {\bf c}$ and $l({\bf a}) = l({\bf c})$.
Conversely, assume ${\bf a} \trianglelefteq {\bf c}$ and
$l({\bf a}) = l({\bf c})$.
This implies ${\bf a}^\circ \trianglelefteq {\bf c}^\circ$ and
$l({\bf a}-1) \geq l({\bf c}-1)$.
We get 
\[ 
\Delta({\bf a}) \subset \overline{\Delta({\bf c})} 
\] 
by applying Lemma \ref{degoftopsocleonemodules}
in case $l({\bf a}-1) > l({\bf c}-1)$, or Proposition \ref{vpiisirreducible}
in case $l({\bf a}-1) = l({\bf c}-1)$.
This finishes the proof.
\end{proof}


\section{Classification of non-regular irreducible 
components}\label{sectionnonregcomp}


The classification of irreducible components of $\V(n,a,b)$ with $a < n$
and $b < n$ is less straightforward than for the case $a = b =n$.
The main reason is that Corollary \ref{stringreduction2}
does not hold in the general case.

A module $M$ is {\it semi-projective} (respectively 
{\it semi-injective}) if it is isomorphic to
\[
\bigoplus_{i=1}^t M(C_i)
\]
where $C_i = x^{a-1}C_i'y^{b-1}$ for some string $C_i'$ and all $i$
(respectively 
$C_i = y^{b-1}C_i'x^{a-1}$ for some string $C_i'$ and all $i$).
The next two statements are clear.

\begin{Lem}\label{semiproj2}
If $M(C)$ is semi-projective and $M(D)$ semi-injective, then
$CxDy$ is a band.
Thus,
\[
M(C) \oplus M(D) \in \overline{{\mathcal F}(CxDy)}.
\]
\end{Lem}

\begin{Lem}\label{nonreg1}
If $M(C)$ is not semi-projective and not semi-injective, then there
exists some $u \in \{ x,y \}$ such that $Cu$ is a band.
Thus, 
\[
M(C) \in \overline{{\mathcal F}(Cu)}.
\] 
\end{Lem}

The next lemma is again a consequence of the construction of
projective covers of string modules.

\begin{Lem}\label{semiprojinj}
Let $M \in \V(n,a,b)$ be a direct sum of $t$ string modules.
If $M$ is semi-projective (respectively semi-injective), then 
$M$ is in some stratum $\sub(L)$ with $\dm \Hom(L,S) = n-t$ (respectively  
$\dm \Hom(L,S) = n+t$).  
\end{Lem}

\begin{Lem}\label{semiproj}
If $M$ is a semi-projective module in $\V(n,a,b)$, then $M$ is
not contained in the closure of the set of regular elements in
$\V(n,a,b)$.
\end{Lem}

\begin{proof}
Let 
\[
M = \bigoplus_{i=1}^t M(C_i)
\]
be semi-projective.
We have $M \in {\mathcal S}(L)$ for some index module $L$.
By Lemma \ref{semiprojinj} we have
\[
\dm \Hom(L,S) = n-t.
\] 
Now assume that ${\mathcal S}(L)$ is contained in the closure of some
stratum ${\mathcal S}(L({\bf a},{\bf b}))$ with $({\bf a},{\bf b})$
regular.
So $L({\bf a},{\bf b}) \leq_{\rm deg} L$.
Since $({\bf a},{\bf b})$ is regular, we get
\[
\dm \Hom({\mathcal S}(L({\bf a},{\bf b}),S) = n.
\]
This is a contradiction because 
the function $\dm \Hom(-,S)$ is upper semicontinuous.
\end{proof}

Lemma \ref{semiproj} enables us to determine when
all irreducible components of $\V(n,a,b)$ are regular, i.e
we can prove Proposition \ref{regulardense}.

\begin{proof}[Proof of Proposition \ref{regulardense}]
If $n \leq a+b-2$ or $n = a+b$, then
there are no semi-projective or semi-injective modules.
So Lemma \ref{nonreg1} implies the result.
For the other direction,
it is sufficient to construct for each $n \geq a+b+1$ and for 
$n = a+b-1$ an $n$-dimensional semi-projective module.
We leave this as an easy exercise to the reader.
Then Lemma \ref{semiproj} yields the result.
\end{proof}

\begin{Lem}\label{reductiontosemiprojectives}
Let $M(C)$ be semi-projective, and let $B$ be a band of the form $x^cy^d$.
Then there exists a semi-projective string module $M(E)$ such that 
\[ 
M(E) \leq_{\rm deg} M(C) \oplus M(B,\lambda). 
\]
\end{Lem}

\begin{proof}
Let $B = x^cy^d$
for some $1 \leq c \leq a-1$ and $1 \leq d \leq b-1$, and let
\[ 
C = x^{c_1}y^{d_1} \cdots x^{c_t}y^{d_t} 
\]
where $1 \leq c_i \leq a-1$ and $1 \leq d_i \leq b-1$ for all $i$, 
$c_1 = a-1$ and $d_t = b-1$.
Note that $M(C)$ is semi-projective.
Let $m$ be the maximal $i$ such that one of the 
following hold:
\begin{enumerate}

\item[(1)] $c_i > c$;

\item[(2)] $c_i = c$ and $d_{i-1} < d$;

\item[(3)] $i = 1$.

\end{enumerate}
First, we assume that there exists some $i \geq m$ such that $d_i < d$.
Note that this implies $i < t$.
Then it follows from the definition of $m$ that $c_{i+1} < c$. 
Define 
\[ 
E = x^{c_1}y^{d_1} \cdots y^{d_i} x^cy^d x^{c_{i+1}} \cdots x^{c_t}
y^{d_t}. 
\]
Now it is easy to construct a short exact sequence
\[ 
0 \longrightarrow M(C) \longrightarrow M(E) \longrightarrow M(B,\lambda) 
\longrightarrow 0. 
\]
This implies $M(E) \leq_{\rm deg} M(C) \oplus M(B,\lambda)$.

Second, we consider the case $d_i \geq d$ for all $i \geq m$.
Let $l$ be maximal such that 
\[ 
C = x^{c_1}y^{d_1} \cdots y^{d_{m-1}} x^{c_m-c} (x^cy^d)^l D 
\]
for some string $D$.
Define 
\[ 
E = x^{c_1}y^{d_1} \cdots y^{d_{m-1}} x^{c_m-c} (x^cy^d)^{l+1} D. 
\]
Again, one can construct a short exact sequence 
\[ 
0 \longrightarrow M(C) \longrightarrow M(E) \longrightarrow M(B,\lambda) 
\longrightarrow 0 
\]
which implies $M(E) \leq_{\rm deg} M(C) \oplus M(B,\lambda)$. 
This finishes the proof.
\end{proof}

Let ${\mathcal P}_n$ (respectively ${\mathcal I}_n$) be the set of all 
semi-projective (respectively semi-injective) modules in $\V(n,a,b)$.
Observe that ${\mathcal P}_n$ and ${\mathcal I}_n$ contain only finitely many
isomorphism classes of modules.
The next corollary follows from Proposition \ref{reductiontodiamonds},
Lemmas \ref{semiproj2}, \ref{nonreg1} and \ref{reductiontosemiprojectives}.

\begin{Cor}\label{nonreg5}
Each non-regular irreducible component of $\V(n,a,b)$ 
contains a dense orbit $\orb$ with 
$\orb \subset {\mathcal P}_n \cup {\mathcal I}_n$.
\end{Cor}

Note that the duality 
${\rm D} = \Hom_\field(-,\field)$ 
induces an isomorphism  
\[
\theta: \V(n,a,b) \longrightarrow \V(n,a,b)
\]
\[
(A,B) \mapsto (A^t,B^t)
\]
where $M^t$ denotes the transpose of a matrix $M$.
For example, if $(A,B)$ is isomorphic to $M(xxy)$, then
$\theta(A,B)$ is isomorphic to ${\rm D}M(xxy) = M(yxx)$.
The restriction of $\theta$ to ${\mathcal P}_n$ yields
an isomorphism ${\mathcal P}_n \to {\mathcal I}_n$.

\begin{Lem}\label{semiprojstratum}
Let ${\mathcal S}(L)$ be a stratum containing a semi-projective module,
and let ${\mathcal S}(M)$ be a stratum containing a semi-injective module.
Then
\[
\sub(L) \not\subseteq \overline{\sub(M)} \text{ and }
\sub(M) \not\subseteq \overline{\sub(L)}. 
\]
\end{Lem}

\begin{proof}
By Lemma \ref{semiprojinj} we get 
\[
\dm \Hom(L,S) = n-s
\] 
and 
\[
\dm \Hom(M,S) = n+t
\]
for some $s,t \geq 1$.
This implies $M \not\leq_{\rm deg} L$.
Thus by Theorem \ref{richmond},(2) the stratum $\sub(L)$ cannot be 
contained in the closure of $\sub(M)$.
Next, assume that $\sub(M)$ is contained in the closure of
$\sub(L)$.
This implies that $\theta(\sub(M))$ is contained in the closure of
$\theta(\sub(L))$ with
$\theta(\sub(M))$ containing a semi-projective and $\theta(\sub(L))$
containing a semi-injective module.
But this is a contradiction to the first part of the proof.
\end{proof}

Up to now, we established the following:
To classify all non-regular irreducible components of 
$\V(n,a,b)$, it is sufficient to decide which orbits in ${\mathcal P}_n$ 
are open.

Let $X$ be indecomposable and semi-projective, and assume that
$X$ is contained in a stratum $\sub(L)$.
We want to determine when $\orb(X)$ is open.
We can assume that $L$ is flip minimal, otherwise we could use
flip degenerations and Theorem \ref{richmond} to show that 
$\sub(L)$ and in particular $X$ is contained in the closure of
some other stratum $\sub(M)$ with $M$ being flip minimal.

Let $({\bf a},{\bf b}) \in {\mathcal P}(n,a,b)$
such that the following hold:
\begin{itemize}
\item $|a \in {\bf a}|, |b \in {\bf b}| \geq 1$;

\item $l({\bf a})+l({\bf b}) = n+1$;

\item $l({\bf a}-1) = l({\bf b}-1)$.

\end{itemize}
Let ${\bf a}-1 = (c_1, \cdots, c_t)$, 
${\bf b}-1 = (d_1, \cdots, d_t)$,
and define
\[ 
P({\bf a},{\bf b}) = M(x^{c_1}y^{d_t}x^{c_2}y^{d_{t-1}} \cdots 
y^{d_2}x^{c_t}y^{d_1}) 
\]
and
\[ 
L({\bf a},{\bf b}) = \Lambda^{n-t} \oplus 
\bigoplus_{i=2}^t M(x^{a-c_i-1}y^{b-d_{t-i+2}-1}). 
\]
Note that $L({\bf a},{\bf b})$ is an index module in 
$\imod_\Lambda(n)$, and
$P({\bf a},{\bf b})$ is semi-projective and contained 
in the stratum $\sub(L({\bf a},{\bf b}))$.
Observe also that $P({\bf a},{\bf b}) \in \Delta({\bf a},{\bf b})$. 
The index module $L({\bf a},{\bf b})$ is flip minimal.
Furthermore, each flip minimal index module $L$ with
$\sub(L)$ containing an indecomposable semi-projective module
is obtained in this way.

\begin{Lem}\label{6.12}
Under the above assumptions,
the orbit $\orb(P({\bf a},{\bf b}))$ is dense
in $\sub(L({\bf a},{\bf b}))$.
\end{Lem}

\begin{proof}
Using the dimension formula in Theorem \ref{richmond},(1) and
Theorem \ref{graphmaps}, 
a straightforward calculation shows that 
\[ 
\dm \orb(P({\bf a},{\bf b})) = 
\dm \sub(L({\bf a},{\bf b})). 
\]
Thus $\orb(L({\bf a},{\bf b}))$ must be dense in
the stratum $\sub(L({\bf a},{\bf b}))$.
\end{proof}

As a consequence of the above results we get the following:

\begin{Lem}\label{opencrit}
The orbit $\orb(P({\bf a},{\bf b}))$ is open if and only if there is no
module $P({\bf c},{\bf d})$ with
$P({\bf c},{\bf d}) <_{\rm deg} P({\bf a},{\bf b})$.
\end{Lem}

\begin{Lem}\label{6.11}
Let $({\bf a},{\bf b}), ({\bf c},{\bf d}) \in  
{\mathcal P}(n,a,b)$
such that 
\begin{itemize}

\item $|a \in {\bf a}|, |b \in {\bf b}|,
|a \in {\bf c}|$, $|b \in {\bf d}| \geq 1$;

\item $l({\bf a}) + l({\bf b}) = l({\bf c}) + l({\bf d}) = n+1$;

\item $l({\bf a}-1) = l({\bf b}-1)$ and
$l({\bf c}-1) = l({\bf d}-1)$.

\end{itemize}
Then the following hold:
\begin{itemize}

\item[(1)]
If $P({\bf c},{\bf d}) \leq_{\rm deg} P({\bf a},{\bf b})$, 
then $({\bf a},{\bf b}) \trianglelefteq ({\bf c},{\bf d})$;

\item[(2)]
If $({\bf a},{\bf b}) \trianglelefteq ({\bf c},{\bf d})$
and $l({\bf a}-1) = l({\bf c}-1)$, then 
$P({\bf c},{\bf d}) \leq_{\rm deg} P({\bf a},{\bf b})$. 

\end{itemize}
\end{Lem}

\begin{proof}
The first part of the lemma is a direct consequence of Lemma 
\ref{strataintersection}.
Next, one easily checks that the conditions
$({\bf a},{\bf b}) \trianglelefteq ({\bf c},{\bf d})$ and
$l({\bf a}-1) = l({\bf c}-1)$ allow a sequence of box move degenerations 
\[ 
L({\bf c},{\bf d}) = L_1 \leq_{\rm deg} L_2 \leq_{\rm deg}
\cdots \leq_{\rm deg} L_t = L({\bf a},{\bf b}) 
\]
such that
$\dm \Hom(L_i,\Lambda) = \dm \Hom(L_1,\Lambda)$
for all $i$.
As before we use Theorem \ref{richmond},(3) and get 
\[ 
\sub(L({\bf a},{\bf b})) \subseteq 
\overline{\sub(L({\bf c},{\bf d}))}. 
\]
Since $P({\bf a},{\bf b})$ and 
$P({\bf c},{\bf d})$ are dense in $\sub(L({\bf a},{\bf b}))$ and
$\sub(L({\bf c},{\bf d}))$, respectively, this implies 
$P({\bf c},{\bf d}) \leq_{\rm deg} P({\bf a},{\bf b})$. 
This finishes the proof.
\end{proof}

\begin{Thm}[Classification of open orbits]\label{openorbits}
Let $X$ be an indecomposable $\Lambda$-module.
Then $\orb(X)$ is open in $\V(n,a,b)$ if and only if $X$ is isomorphic to 
$M(C)$ or ${\rm D}M(C)$ where $C$ is of one of 
the following forms:
\begin{enumerate}

\item[{\rm (1)}]
\[
C = (x^{a-1}y)^r(x^{a-1}y^{b-1})^s(xy^{b-1})^t 
\]
where $r,s,t \geq 0$, $r+s \geq 1$ and $s+t \geq 1$;

\item[{\rm (2)}]
\[ 
C = (x^{a-1}y)^r (x^{a-1}y^i)^\alpha (x^{a-1}y^{b-1})^s 
(x^jy^{b-1})^\beta (xy^{b-1})^t 
\]
where $r,s,t \geq 0$, $2 \leq i \leq b-2$, $2 \leq j \leq a-2$,
$0 \leq \alpha, \beta \leq 1$, $\alpha + \beta \geq 1$, 
$r+\alpha+s \geq 1$ and $s+\beta+t \geq 1$;

\item[{\rm (3)}] 
\[ 
C = (x^{a-1}y)^r x^iy^j (xy^{b-1})^t 
\]
where $r,t \geq 1$, $1 \leq i \leq a-2$ and $1 \leq j \leq b-2$.

\end{enumerate}
The open orbits in $\V(n,a,b)$ are exactly the orbits
of the form
\[ 
\orb \left( \bigoplus\limits_{i \in I} M(C_i) \right) 
\] 
with $\orb(M(C_i))$ open and $\Ext^1(M(C_i),M(C_j)) = 0$ 
for all $i \not= j$ in $I$.
\end{Thm}

If a string $C$ belongs to one of the sets (1), (2) or (3) as 
defined in the theorem, then we say that $C$ is of {\it type} 
(1), (2) or (3), respectively.

\begin{proof}
We classify the open orbits $\orb(X)$ with $X$ indecomposable.
By Lemma \ref{nonreg1} we know that $X$ 
has to be semi-projective or semi-injective.
By duality, we can assume without loss of generality that $X$ is 
semi-projective.
As a consequence of Lemma \ref{6.11}, we can assume that 
$X = M(C) = P({\bf a},{\bf b})$ such that the
following hold:
\begin{itemize}

\item ${\bf a}$ has at most one entry different from $1,2$ and $a$;

\item ${\bf b}$ has at most one entry different from $1,2$ and $b$.

\end{itemize}
Now we proceed similar to the proof of Theorem
\ref{theorem1}.
We can assume that $l({\bf a}) + l({\bf b}) = n+1$,
\[
{\bf a}-1 = ((a-1)^{p-r-1},a-v-1,1^r)
\] 
and
\[
{\bf b}-1 = ((b-1)^{p-s-1},b-w-1,1^s)
\]
where 
$0 \leq v \leq a-2$, $0 \leq w \leq b-2$, $0 \leq r,s \leq p-1$,
$v=0 \Rightarrow r=0$ and $w=0 \Rightarrow s=0$.
Then by using Theorem \ref{graphmaps}, we get  
\begin{multline*}
\dm \orb(P({\bf a},{\bf b})) =
n^2 - p^2-p-1 - (a-v-2)(p-r)^2\\ 
- (b-w-2)(p-s)^2
- v(p-r-1)^2 - w(p-s-1)^2.
\end{multline*}
By Lemma \ref{opencrit}
the orbit of $P({\bf a},{\bf b})$ is open if and only if there is no
$P({\bf c},{\bf d})$ with $
P({\bf c},{\bf d}) <_{\rm deg} P({\bf a},{\bf b})$.

If $r+s+1 \leq p$, then 
\[
\dm \orb(P({\bf a},{\bf b})) 
\geq \dm \orb(P({\bf c},{\bf d}))
\]
for all 
$P({\bf c},{\bf d})$ with $({\bf a},{\bf b}) \trianglelefteq
({\bf c},{\bf d})$.
This follows from the above dimension formula.
So by Lemma \ref{6.11} the orbit 
$\orb(P({\bf a},{\bf b}))$ must be open.
Observe that $C$ is of type (1),(2) or (3) if and only if 
$r+s+1 \leq p$.

Next, assume that $r+s+1 > p$.
By the definition of $r$ and $s$, 
it follows that $a,b \geq 3$ in this case.
Then $C$ is of the form
\[ 
(x^{a-1}y)^k (x^iy)(xy)^l(xy^j)(xy^{b-1})^m 
\]
where $k,m \geq 1$, $l \geq 0$, $1 \leq i \leq a-2$ and 
$1 \leq j \leq b-2$.
If $l=0$, then define
\[ 
E = (x^{a-1}y)^kx^{i+1}y^{j+1}(xy^{b-1})^m. 
\]
Otherwise, let 
\[ 
E = (x^{a-1}y)^k(x^{i+1}yy)(xy)^{l-1}(xy^j)(xy^{b-1})^m. 
\]
In both cases, we get $M(E) <_{\rm deg} M(C)$.
This is proved by constructing a Riedtmann sequence
\[ 
0 \longrightarrow M(C) \longrightarrow M(E) \oplus M(xy,1) \longrightarrow
M(xy,1) \longrightarrow 0. 
\]
Thus, $\orb(P({\bf a},{\bf b}))$ cannot be open.
This finishes the classification of indecomposable
$\Lambda$-modules whose orbit is open.
The rest of the theorem follows from \cite[Theorem 3]{Z}.
\end{proof}

For modules $X$ and $Y$ let $\underline{\Hom}(X,Y)$ be
the space $\Hom(X,Y)$ modulo the homomorphisms factoring 
through a projective module. 
By $\tau$ we denote the Auslander-Reiten translation.
For indecomposable modules $X$ and $Y$ we have the Auslander-Reiten
formula
\[ 
\Ext^1(X,Y) \cong {\rm D\underline{Hom}}(\tau^{-1}Y,X). 
\]
For the basics of Auslander-Reiten theory we refer to 
\cite{ARS} or \cite{Ri}.
If $M(C)$ is a semi-projective string module, then define 
\[ 
\tau^{-1}C = x^{a-1}yCxy^{b-1}. 
\]
Note that $M(\tau^{-1}C)$ is also semi-projective.
It is proved in \cite{BR} that 
\[ 
\tau^{-1}M(C) = M(\tau^{-1}C). 
\]
The next proposition is an application of the Auslander-Reiten formula
and Theorem \ref{graphmaps}.

\begin{Prop}\label{extzerosemipro}
If $M(C)$ and $M(D)$ are semi-projective string modules, then
the following are equivalent:
\begin{itemize}

\item[(1)]
$\Ext^1(M(C),M(D)) = 0$;

\item[(2)] 
Each map $f_a$ with $a \in {\mathcal A}(\tau^{-1}D,C)$ 
factors through $M(x^{a-1}y^{b-1})$.

\end{itemize}
\end{Prop}

For deciding whether a graph map factors through another string module, 
one uses the multiplicative behaviour of graph maps.
Using this proposition and the previous theorem, 
it is now easy to compute the semi-projective modules whose orbit is 
open.
Using duality, we get all open orbits. 
This completes the classification of irreducible components of the variety
$\V(n,a,b)$.

\begin{Cor}\label{extzeromodules}
For an indecomposable
$\Lambda$-module $X$ the following statements 
are equivalent:
\begin{itemize}

\item[(1)] $\Ext^1(X,X) = 0$;

\item[(2)]
$X$ is isomorphic to a string module $M(C)$ or 
${\rm D}M(C)$ with 
\[ 
C = (x^{a-1}y)^r (x^{a-1}y^{b-1})^s (xy^{b-1})^t 
\]
where $r,s,t \geq 0$, $r+s \geq 1$ and $s+t \geq 1$.

\end{itemize}
\end{Cor}


\section{Remarks and examples}\label{examples}


We list all irreducible components of $\V(n,3,3)$ for $n \leq 12$.
First, let us give the list of all regular irreducible components and their
dimensions.

For each regular $({\bf a},{\bf b})$ we constructed a family
${\mathcal F}({\bf a},{\bf b})$ of modules which is dense in
$\Delta({\bf a},{\bf b})$, see Proposition \ref{reductiontodiamonds}.
Recall that these families are of the form 
${\mathcal F}((B_1,p_1), \cdots, (B_m,p_m))$.

In Figure \ref{fig6}
we display the data $(B_1,p_1),\cdots,(B_m,p_m)$ 
in case the closure of the corresponding family is an irreducible component.
If $p_i = 1$, then we just write $B_i$ instead of $(B_i,p_i)$.
\begin{figure}[ht]
\begin{tabular}{rl|rl|rl|rl|rl}
2     & & 3      & & 4        &  & 5         &  & 6             &\\
\hline
\hline
$xy$&3& $xxy$&7& $xxyy$ &13& $xxy,xy$&20& $(xxy,2)$ &28\\
      & & $xyy$&7&          &  & $xyy,xy$&20&  $(xyy,2)$   &28\\
      & &        & &          &  &           &  & $xxy,xyy$ & 30  
\end{tabular}

\vspace{0.7cm}

\begin{tabular}{rl|rl|rl}
7           &  & 8               &  & 9                &\\
\hline
\hline
$xxyy,xxy$&40  & $(xxy,2),xy$&51   & $(xxy,3)$     & 63 \\
$xxyy,xyy$&40  & $(xyy,2),xy$&51   & $(xyy,3)$     & 63 \\
            &  & $(xxyy,2)$  &52   & $(xxy,2),xyy$ & 67 \\
            &  & $xxy,xyy,xy$  &53 & $(xyy,2),xxy$ & 67 
\end{tabular}

\vspace{0.7cm}

\begin{tabular}{rl|rl|rl}
10                &  & 11                  &   & 12                     &\\
\hline
\hline
 $(xxy,2),xxyy$ &81   & $(xxy,3),xy$     &96  & $(xxy,4)$          &112 \\
 $(xyy,2),xxyy$ &81   & $(xyy,3),xy$     &96  & $(xyy,4)$          &112 \\
 $xxyy,xxy,xyy$ &83   & $(xxyy,2),xxy$   &99  & $(xxyy,3)$         &117 \\ 
                   &  & $(xxyy,2),xyy$   &99  & $(xxy,3),xyy$      &118 \\
                   &  & $(xxy,2),xyy,xy$ &100 & $(xyy,3),xxy$      &118 \\
                   &  & $(xyy,2),xxy,xy$ &100 & $(xxy,2),(xyy,2)$  &120  
\end{tabular}
\caption{The regular components of $\V(n,3,3)$ for $n \leq 12$}\label{fig6}
\end{figure}

\vspace{0.4cm}
In Figure \ref{fig7} we give a list of all open
orbits and their dimensions.
Recall that the closures of the open orbits are exactly the non-regular 
irreducible components.
Remember also that the open orbits are orbits of certain semi-projective 
or semi-injective modules.
For the sake of brevity we list only the strings $C_i$ occurring in their 
direct sum decomposition.
For example $xxyy \oplus xxyy$ encodes the module $M(xxyy) \oplus M(xxyy)$.
We only list the semi-projective modules whose orbits are open. 
Thus one has to add the same number of semi-injective modules to get 
all open orbits.
Recall that there are no open orbits for $n \leq a+b-2 = 4$ and
$n = a+b = 6$.
\begin{figure}[ht]
\begin{tabular}{rl|rl|rl|rl} 
5     &  & 7       &  & 8        &  & 9           &\\
\hline
\hline
$xxyy$&20& $xxyxyy$&40& $xxyyxyy$&52& $(xxyy)^2$  &66\\
      &  &         &  & $xxyxxyy$&52& $xxyxyxyy$  &66
\end{tabular}

\vspace{0.7cm}

\begin{tabular}{rl|rl|rl} 
10                &  & 11           &   & 12                  &\\
\hline
\hline
$xxyy \oplus xxyy$&80 & $(xxy)^2xxyy$ &98  & $xxyy \oplus xxyxyy$ &117 \\
$(xxy)^2xyy$      &82 & $xxyy(xyy)^2$ &98  & $(xxyy)^2xyy$        &118 \\
$xxy(xyy)^2$      &82 & $xxyxxyyxyy$  &100 & $xxy(xxyy)^2$        &118  \\
                  &   &              &     & $(xxy)^2xyxyy$       &118  \\
                  &   &              &     & $xxyxy(xyy)^2$       &118  
\end{tabular}
\caption{The non-regular components of $\V(n,3,3)$ for $n \leq 12$}
\label{fig7}
\end{figure}

\vspace{0.4cm}
{\bf Remark 1}:
If $\Ext^1(M,M) = 0$ for some $\Lambda$-module $M$, then by Voigt's Lemma
one gets that $\orb(M)$ is open.
The converse does not hold.
The smallest example of this kind occurs for $n=9$:
Let $M = M(xxyxyxyy)$ be in $\V(9,3,3)$.
Then $\orb(M)$ is open but $\Ext^1(M,M) \not= 0$.

\vspace{0.3cm}
{\bf Remark 2}:
Let $a = b = 2$ and $n=3$.
Then $\orb(M(xy))$ and $\orb(M(yx))$ are both open orbits,
since $M(xy)$ is projective and $M(yx)$ is injective.
In particular, $\Delta(2,1)$ and $\Delta((2,1),(2,1))$ are both not
irreducible.

\vspace{0.3cm}
{\bf Remark 3}:
The Gelfand-Ponomarev algebra $\Lambda$ is a string algebra in the 
sense of \cite{BR}.
Similarly to Lemma \ref{stringalgebrasaresubfinite}  
one can show that all string algebras are subfinite,
and their index modules can be classified as in Lemma \ref{indexmodules}.

One should be able to classify the irreducible components
of varieties of modules over many other string algebras in the same fashion
as in this paper.


\end{document}